\documentstyle[twocolumn,eqsecnum,aps]{revtex}

\newcommand{\Ent}[1]{[\mkern - 2.5 mu [#1] \mkern - 2.5 mu ]}

\begin{document}

\draft

\title{Irregular Input Data in Convergence Acceleration and Summation 
Processes: General Considerations and Some Special Gaussian
Hypergeometric Series as Model Problems}

\author{Ernst Joachim Weniger\cite{Internet}}
\address{Institut f{\" u}r Physikalische und Theoretische Chemie \\
Universit{\" a}t Regensburg, D-93040 Regensburg, Federal Republic of
Germany}

\date{Submitted to Computer Physics Communications --
2 March 2000}
\maketitle

\begin{abstract}
Sequence transformations accomplish an acceleration of convergence or a
summation in the case of divergence by detecting and utilizing
regularities of the elements of the sequence to be transformed. For
sufficiently large indices, certain asymptotic regularities normally do
exist, but the leading elements of a sequence may behave quite
irregularly. The Gaussian hypergeometric series ${}_2 F_1 (a, b; c; z)$
is well suited to illuminate problems of that kind. Sequence
transformations perform quite well for most parameters and
arguments. If, however, the third parameter $c$ of a nonterminating
hypergeometric series ${}_2 F_1$ is a negative real number, the terms
initially grow in magnitude like the terms of a mildly divergent
series. The use of the leading terms of such a series as input data
leads to unreliable and even completely nonsensical results. In
contrast, sequence transformations produce good results if the leading
irregular terms are excluded from the transformation process. Similar
problems occur also in perturbation expansions. For example, summation
results for the infinite coupling limit $k_3$ of the sextic anharmonic
oscillator can be improved considerably by excluding the leading terms
from the transformation process. Finally, numerous new recurrence
formulas for the ${}_2 F_1 (a, b; c; z)$ are derived.
\end{abstract}


\section{Introduction}
\label{sec:Intro} 

In mathematics and in the mathematical treatment of scientific problems,
slowly convergent or divergent sequences and series occur abundantly.
Accordingly, many techniques for the acceleration of convergence and the
summation of divergent series have been invented, and some of them are
even older than calculus (see for instance pp.\ 90 - 91 of
\cite{Br91} or p.\ 249 of \cite{Kn64}).

{\em Sequence transformations\/} are principal tools to overcome
convergence problems. Let us assume that $\{ s_n \}_{n=0}^{\infty}$ is a
slowly convergent or divergent sequence, whose elements $s_n$ may for
example be the partial sums of an infinite series:
\begin{equation}
s_n \; = \; \sum_{k=0}^{n} \, a_k \, .
\end{equation}
The basic assumption of all sequence transformations is that a sequence
element $s_n$ can for all indices $n \ge 0$ be partitioned into a
(generalized) limit $s$ and a remainder or truncation error $r_n$
according to
\begin{equation}
s_n \; = \; s \, + \, r_n \, .
\end{equation}
The conventional approach of evaluating an infinite series consists in
adding up so many terms that the remainders $r_n$ ultimately become
negligible. Unfortunately, this is not always feasible because of
obvious practical limitations. Moreover, adding up further terms does
not work in the case of a divergent series since their terms usually
increase in magnitude with increasing index.

Alternatively, one could try to determine approximations to the
remainders $r_n$ and to eliminate them from the sequence elements
$s_n$. At least conceptually, this is what a sequence transformation
tries to accomplish. Thus, the original sequence $\{ s_n
\}_{n=0}^{\infty}$ is transformed into a new sequence $\{ s'_n
\}_{n=0}^{\infty}$ whose elements have the same (generalized) limit $s$ 
but different remainders $r'_n$:
\begin{equation}
s'_n \; = \; s \, + \, r'_n \, .
\end{equation}
The transformation process was successful if the transformed remainders
$r'_n$ have superior numerical properties. For example, in the
literature on extrapolation methods it is said that a sequence
transformation accelerates convergence if the transformed remainders
$r'_n$ vanish more rapidly than the original remainders $r_n$ according
to
\begin{equation}
\lim_{n \to \infty} \, \frac{r'_n}{r_n} \; = \; 
\lim_{n \to \infty} \, \frac{s'_n - s}{s_n - s} \; = \; 0 \, .
\label{CondConvAcc}
\end{equation}
Similarly, a divergent sequence $\{ s_n \}_{n=0}^{\infty}$, whose
remainders $r_n$ do not vanish as $n \to \infty$, is transformed into
convergent sequence $\{ s'_n \}_{n=0}^{\infty}$ if the transformed
remainders $r'_n$ vanish as $n \to \infty$.

During the last years, considerable progress has been reached in this
field, as documented by the large number of recent monographs
\cite{Bre77,Bre78,Bre80,BreRZa91,Wi81,MaSha83,De88,LiLuSh95} and review
articles \cite{Gu89,We89}. Moreover, numerous applications of sequence
transformations have been reported in the literature. For example, the
present author has applied sequence transformations successfully in such
diverse fields as the evaluation of special functions
\cite{We89,We90,WeCi90,We94,We96d,JeMoSoWe99}, the evaluation of 
molecular multicenter integrals of exponentially decaying functions
\cite{GroWeSte86,HoWe95,WeGroSte86,WeSte89,SteWe90}, the summation of 
strongly divergent quantum mechanical perturbation expansions
\cite{We90,WeCiVi91,We92,WeCiVi93,We96a,We96b,CiWeBraSpi96,We96c,We97,%
JeBeWeS02000}, and the extrapolation of crystal orbital and cluster
calculations for oligomers to their infinite chain limits of stereoregular
quasi-onedimensional organic polymers \cite{WeLi90,CioWe93}.

It should be noted that Pad\'{e} approximants \cite{BaGM96}, which in
applied mathematics and in theoretical physics have become the standard
tool to overcome convergence problems with power series, can be
considered to be a special class of sequence transformations since the
partial sums of a power series are transformed into a doubly indexed
sequence of rational functions.

As described above, sequence transformations try at least in principle
to construct approximations to the actual remainders which are then
eliminated from the input data. This is done by detecting and utilizing
regularities in the behavior of the elements of the sequence to be
transformed. For sufficiently large indices $n$, one can expect that
certain asymptotic regularities do exist. However, sequence
transformations are normally used with the intention of avoiding the
asymptotic domain, i.e., the transforms are constructed from the
\emph{leading} elements of the input sequence. Unfortunately, sequence
elements $s_n$ with small indices $n$ often behave irregularly. In such
a case, a straightforward application of a sequence transformation may
be ineffective and even lead to completely nonsensical results. Instead,
one should analyze the behavior of the input data as a function of the
index and exclude highly irregular sequence elements from the
transformation process if necessary. In this way, one has a much better
chance of obtaining good and reliable transformation results. It is the
intention of this article to describe and classify some of the problems,
which can result from irregular input data, and to discuss strategies to
overcome them.

Sequence transformations are needed most in cases in which apart from
the numerical values of a few elements of a slowly convergent or
divergent sequence only very little is known. This is a situation which
is not uncommon in scientific applications as for example the summation
of strongly divergent perturbation expansions as they occur quantum
mechanics or in quantum field theory. Thus, it would in principle be
desirable to discuss problems with irregular input data also via
examples of numerically determined input data. However, the lack of
detailed knowledge about the behavior of the elements of such a sequence
makes it hard to fully understand the numerical problems as well as to
develop strategies to overcome them. Consequently, complications of that
kind are discussed in this article predominantly via suitable
mathematical model problems.

In Section \ref{sec:Paths}, some formal aspect of sequence
transformations are discussed, in particular the concept of a {\em
path\/} in the table of a transformation. In this way, different
approaches for the determination of the approximations to the
(generalized) limit of the input sequence can be classified and
formalized. Moreover, the concept of a path helps to understand the
impact of irregular input data on the performance of sequence
transformations.

At least some of the problems mentioned above can be illuminated by
considering the Gaussian hypergeometric function ${}_2 F_1 (a, b; c; z)$
which is defined by a power series that converges in the interior of the
unit circle. This function does not only depend on an argument $z$ but
also on three essentially arbitrary parameters $a$, $b$, and $c$. As
discussed in Section \ref{sec:2F1}, the convergence of the
hypergeometric series can for most values of the parameters $a$, $b$,
and $c$ be accelerated quite effectively by a variety of different
sequence transformations. Moreover, it is in this way frequently
possible to associate a finite value to a hypergeometric series even if
its argument does not lie in the interior of the unit circle.

However, as discussed in Section \ref{sec:Ne3Pa2F1}, the situation
changes dramatically if the third parameter $c$ of the hypergeometric
series ${}_2 F_1$ is a negative real number. Then, the terms of this
series first increase with increasing summation index even for $\vert z
\vert < 1$ and produce partial sums which look like the elements of a 
mildly divergent sequence. Only for sufficiently large indices, the
terms decrease in magnitude and ultimately produce a convergent
result. Accordingly, the leading partial sums of such a hypergeometric
series display a highly irregular behavior, and they should not be used
as input data for a sequence transformation. If the leading irregular
coefficients are skipped and only regular coefficients with higher
indices are used as input data, then sequence transformations are again
able to produce good and reliable results.

Section \ref{sec:SumConclu} contains a summary. In Appendix
\ref{App_A}, the properties of the sequence transformations, which are
used in this article, are discussed. Problems with irregular input data
occur also quite frequently in the mathematical treatment of scientific
problems. In Appendix \ref{App_B}, it is shown that the convergence of
extensive summation calculations for the so-called infinite coupling
limit $k_3$ of the sextic anharmonic oscillator can be improved
considerably by excluding the leading irregular coefficients of the
divergent perturbation series from the transformation process. Finally,
numerous new recurrence formulas for the hypergeometric function ${}_2
F_1 (a, b; c; z)$ are derived in Appendix \ref{App_C}.

\section{Order-Constant and Index-Constant Paths}
\label{sec:Paths} 

In this Section, some aspects of sequence transformations are discussed
which admittedly look very formal. Nevertheless, they should not be
ignored since they may be very consequential in practical applications.

Obviously, a computational algorithm can only involve a finite number of
arithmetic operations. Consequently, a sequence transformation ${\cal
T}$ can only use finite subsets of the original sequence $\{ s_n
\}_{n=0}^{\infty}$ for the computation of new sequence elements
$s^{\prime}_m$. In addition, these finite subsets normally consist of
consecutive elements. Accordingly, only subsets of the type $\{ s_n,
s_{n+1}, \ldots, s_{n+l} \}$ will be considered in this article.

All the commonly used sequence transformations ${\cal T}$ can be
represented by infinite sets of doubly indexed quantities $T_{k}^{(n)}$
with $k, n \ge 0$ that can be displayed in a two-dimensional array which
is called the {\em table\/} of ${\cal T}$.

Here, the convention is used that the superscript $n$ always indicates
the minimal index occurring in the finite subset of sequence elements
used for the computation of a given $T_k^{(n)}$. The subscript $k$ --
usually called the {\em order\/} of the transformation -- is a measure
for the complexity of the transformation process which yields
$T_k^{(n)}$.

The elements $T_k^{(n)}$ of the table of ${\cal T}$ are gauged in such a
way that $T_0^{(n)}$ corresponds to an untransformed sequence element,
\begin{equation}
T_0^{(n)} \; = \; s_n \, .
\end{equation}
An increasing value of $k$ implies that the complexity of the
transformation process increases. Moreover, $l = l(k)$ also
increases. This means that for every $k, n \ge 0$ the sequence
transformation ${\cal T}$ produces a new transform according to
\begin{equation}
 T_k^{(n)} \; = \;
{\cal T} \bigl(s_n, s_{n+1}, \ldots , s_{n+l(k)} \bigr) \, .
\end{equation}
The exact relationship, which connects $k$ and $l$, is specific for a
given sequence transformation ${\cal T}$.

Let us assume that a sequence transformation ${\cal T}$ should be used
to speed up the convergence of some sequence $\{ s_n \}_{n=0}^{\infty}$
to its limit $s = s_{\infty}$. One can try to obtain a better
approximation to $s$ by proceeding on an in principle unlimited variety
of different {\em paths\/} in the table of ${\cal T}$. Two extreme types
of paths -- and also those which are predominantly used in practical
applications -- are {\em order-constant\/} paths
\begin{equation}
T_{k}^{(n)}, T_{k}^{(n+1)}, T_{k}^{(n+2)}, \ldots
\label{OrdConstPath}
\end{equation}
with fixed transformation order $k$ and $n \to \infty$, and {\em
index-constant\/} paths
\begin{equation}
T_{k}^{(n)}, T_{k+1}^{(n)}, T_{k+2}^{(n)}, \ldots
\label{IndConstPath}
\end{equation}
with fixed minimal index $n$ and $k \to \infty$.

Order-constant and index-constant paths differ significantly. It is not
even {\em a priori\/} clear that these two types of paths lead to the
same limit in the case of an arbitrary sequence $\{ s_n
\}_{n=0}^{\infty}$. However, for the sake of simplicity this potential
complication will be ignored here, and we shall always tacitly assume
that order-constant and index-constant paths lead to the same limit.

In the case of an order-constant path, a {\em fixed\/} number of $l+1$
sequence elements $\{ s_n, s_{n+1}, \ldots s_{n+l} \}$ is used for the
computation of $T_k^{(n)}$, and the starting index $n$ of this string of
fixed length is increased successively until either convergence is
achieved or the number of available elements of the input sequence is
exhausted.

In the case of an index-constant path, the starting index $n$ is kept
fixed at a low value (usually $n = 0$ or $n = 1$) and the transformation
order $k$ is increased and with it the number of elements contained in
the subset $\{ s_n, s_{n+1}, \ldots s_{n+l(k)} \}$. Thus, on an
index-constant path it is always tried to compute from a given set of
input data that element $T_k^{(n)}$ which has the highest possible
transformation order $k$.

In order to clarify the differences between order-constant and
index-constant paths, let us consider the probably best known sequence
transformation, Wynn's epsilon algorithm \cite{Wy56}:
\begin{mathletters}
\label{eps_al}
\begin{eqnarray}
\epsilon_{-1}^{(n)} & = & 0 \, , \qquad
\epsilon_0^{(n)} \; = \; s_n \, , \\
\epsilon_{k+1}^{(n)} & = & \epsilon_{k-1}^{(n+1)} \, + \,
1 / [\epsilon_{k}^{(n+1)} - \epsilon_{k}^{(n)} ] \, .
\end{eqnarray}
\end{mathletters}%
\noindent
Wynn \cite{Wy56} showed that if the input data $s_n$ for the epsilon
algorithm are the partial sums
\begin{equation}
f_n (z) \; = \; \sum_{\nu=0}^{n} \gamma_{\nu} z^{\nu}
\label{PowSerPS}
\end{equation}
of a (formal) power series for some function $f (z)$, then the elements
$\epsilon_{2 k}^{(n)}$ with {\em even\/} subscripts are Pad\'{e}
approximants to $f$ according to
\begin{equation}
\epsilon_{2 k}^{(n)} \; = \; [ n + k / k ] \, .
\label{Eps_Pade}
\end{equation}
Here, the notation of the monograph by Baker and Graves-Morris
\cite{BaGM96} is used, i.e., a Pad\'{e} approximant $[ l / m ]$
corresponds to the ratio of two polynomials $P_l (z)$ and $Q_m (z)$,
which are of degrees $l$ and $m$, respectively, in $z$. In contrast, the
elements $\epsilon_{2 k + 1}^{(n)}$ with {\em odd\/} subscripts are only
auxiliary quantities which diverge if the whole process converges.

It follows from (\ref{Eps_Pade}) that the epsilon algorithm
(\ref{eps_al}) effects the following transformation of the partial sums
(\ref{PowSerPS}) to Pad\'{e} approximants:
\begin{equation}
\bigl\{ f_{n} (z), f_{n+1} (z), \ldots, f_{n+2k} (z) \bigr\}
\; \longrightarrow \; [ n + k / k ] \, .
\end{equation}
Thus, if we use a window consisting of $2k + 1$ partial sums $f_{n+j}
(z)$ with $0 \le j \le 2k$ on an order-constant path and increase the
minimal index $n$ successively, the epsilon algorithm produces the
following sequence of Pad\'{e} approximants:
\begin{equation}
[ n + k / k ], [ n + k + 1 / k ], \ldots, [ n + k + m / k ],
\ldots \, .
\end{equation}
Only $2k + 1$ partial sums are used for the computation of the Pad\'{e}
approximants, although many more are known. Obviously, the available
information is not exploited optimally on such an order-constant path.

Moreover, the degree of the numerator polynomial of a Pad\'{e}
approximant $[ n + k + m / k ]$ increases with increasing $m \ge 0$,
whereas the degree of the denominator polynomial remains fixed. Thus,
these Pad\'{e} look unbalanced. Instead, it seems to be much more
natural to use {\em diagonal\/} Pad\'{e} approximants, i.e., Pad\'{e}
approximants with numerator and denominator polynomials of equal degree,
or -- if this is not possible -- to use Pad\'{e} approximants with
degrees of the numerator and denominator polynomials that differ as
little as possible.

This approach has in principle many theoretical as well as practical
advantages. For example, Wynn could show that if the partial sums $f_0
(z)$, $f_1 (z)$, $\cdots$, $f_{2 n} (z)$ of a Stieltjes series are used
for the computation of Pad\'{e} approximants, then the diagonal
approximant $[ n / n ]$ provides the most accurate approximation to the
corresponding Stieltjes function $f (z)$, and if the partial sums $f_0
(z)$, $f_1 (z)$, $\cdots$, $f_{2 n + 1} (z)$ are used for the
computation of Pad\'{e} approximants, then either $[ n + 1 / n ]$ or $[
n / n + 1 ]$ provides the most accurate approximation (Theorem 5 of
\cite{Wy68}). A detailed discussion of Stieltjes series and their
special role in the theory of Pad\'{e} approximants can for instance be
found in Section 5 of the monograph by Baker and Graves-Morris
\cite{BaGM96}.

Thus, it is apparently an obvious idea to try to use either diagonal
Pad\'{e} approximants or their closest neighbors whenever possible. Let
us assume that the partial sums $f_0 (z)$, $f_1 (z)$, $\ldots$, $f_{m}
(z)$ are known. If $m$ is even or odd, $m = 2 \mu$ or $m = 2 \mu + 1$,
respectively, the elements of the epsilon table with the highest
possible transformation orders are given by the transformations
\begin{eqnarray}
\bigl\{ f_0 (z), f_1 (z), \ldots, f_{2\mu} (z) \bigr\}
& \longrightarrow & \epsilon_{2 \mu}^{(0)} = [ \mu / \mu ] \, ,
\\
\bigl\{ f_1 (z), f_2 (z), \ldots, f_{2\mu+1} (z) \bigr\}
& \longrightarrow & \epsilon_{2 \mu}^{(1)} = [ \mu + 1 / \mu ]
\, .
\end{eqnarray}
With the help of the notation $\Ent x$ for the integral part of $x$,
which is the largest integer $\nu$ satisfying $\nu \le x$, these two
relationships can be expressed by a single equation (Eq.\ (4.3-6) of 
\cite{We89}):
\begin{eqnarray}
\lefteqn{\left\{ f_{m - 2 \Ent {m/2}} (z),
f_{m - 2 \Ent {m/2} + 1} (z), \ldots , f_m (z) \right\}}
\qquad \nonumber \\
\qquad & \longrightarrow &
\epsilon_{2 \Ent {m/2}}^{(m - 2 \Ent {m/2})} \; = \;
\bigl[ m - \Ent {m/2} / \Ent {m/2} \bigr] \, .
\label{EpsApprLim}
\end{eqnarray}
For $m = 0, 1, 2, \ldots$, these transformations correspond to the
following staircase sequence in the Pad\'{e} table (Eq.\ (4.3-7) of
\cite{We89}):
\begin{eqnarray}
\lefteqn{[0/0], [1/0], [1/1], \ldots} \nonumber \\
& & \quad \ldots, [\nu / \nu], [\nu + 1/ \nu], [\nu +1/ \nu +1], \ldots \, .
\end{eqnarray}
This staircase sequence exploits the available information optimally if
the partial sums $f_m (z)$ with $m \ge 0$ are computed successively and
if after the computation of each new partial sum the element of the
epsilon table with the highest possible {\em even\/} transformation
order is computed. Moreover, the Pad\'e approximants obtained in this
way look balanced since the degrees of their numerator and denominator
polynomials differ as little as possible.

The example of Wynn's epsilon algorithm strongly indicates that
index-constant paths are at least in principle computationally more
efficient than order-constant paths since they exploit the available
information optimally. This is in general also true for all other
sequence transformations considered in this article.

Another serious disadvantage of order-constant paths is that they cannot
be used for the summation of divergent sequences and series since
increasing $n$ in the set $\{ s_n, s_{n+1}, \ldots, s_{n+l} \}$ of input
data normally only increases divergence.

In view of the examples given above, it is apparently an obvious idea to
use exclusively index-constant paths, and preferably those which start
at a very low index $n$, for instance at $n = 0$ or $n = 1$. This is
certainly a good idea if all elements of the input sequence contain
roughly the same amount of useful information. If, however, the leading
terms of the sequence to be transformed behave irregularly, they cannot
contribute useful information, or -- to make things worse -- they
contribute \emph{wrong} information. In such a case it is usually
necessary to exclude the leading elements of the input sequence from the
transformation process. Thus, it is preferable to use either an
order-constant path or an index-constant path with a sufficiently large
starting index $n$. The use of an order-constant path has the additional
advantage that the diminishing influence of irregular input data with
small indices $n$ should become obvious from the transformation results
as $n$ increases.

Finally, an important theoretical advantage of order-constant paths
should also be mentioned. Normally, it is much easier to perform a
theoretical convergence analysis of a sequence transformation on an
order-constant path than on an index-constant path. As the starting
index $n$ of the string $\{ s_n, s_{n+1}, \ldots, s_{n+l} \}$ of input
data becomes large, asymptotic approximations to the sequence elements
$s_n$ can be used. Often, this greatly simplifies a theoretical
analysis. In the case of index-constant paths, such a simplification is
not possible because not all input data have large indices.
Accordingly, theoretical convergence properties of sequence
transformations are studied almost exclusively on order-constant
paths. Notable exceptions are two articles by Sidi \cite{Si79,Si80}
where the convergence properties of sequence transformations on both
order-constant and index-constant paths are analyzed.

\section{The Gaussian Hypergeometric Function}
\label{sec:2F1} 

The Gaussian hypergeometric function ${}_2 F_1 (a, b; c; z)$ is one of
the most important special functions of mathematical physics, and its
properties are discussed in numerous books, for example in those by
Abramowitz and Stegun \cite{AbSte72}, Erd\'{e}lyi, Magnus,
Oberhettinger, and Tricomi \cite{ErMaObTri53}, Magnus, Oberhettinger,
and Soni \cite{MaObSo66}, Seaborn \cite{Sea91}, Slater \cite{Sla66},
Spanier and Oldham \cite{SpaOld87}, Temme \cite{Tem96}, and Wang and Guo
\cite{WaGu89}. 

The hypergeometric function ${}_2 F_1 (a, b; c; z)$ is defined via the
corresponding hypergeometric series (p.\ 37 of \cite{MaObSo66})
\begin{equation}
{}_2 F_1 (a, b; c; z) \; = \; \sum_{m=0}^{\infty} \,
\frac {(a)_m (b)_m} {(c)_m m!} \, z^m \, ,
\label{Ser_2f1}
\end{equation}
where $(a)_m = \Gamma(a+m)/\Gamma(a)$ is a Pochhammer symbol (see for
example p.\ 3 of \cite{MaObSo66}). The series (\ref{Ser_2f1}) terminates
after a finite number of terms if either $a$ or $b$ is a negative
integer. Otherwise, it converges in the interior of the unit circle,
i.e., for $\vert z \vert < 1$, and it diverges for $\vert z \vert >
1$. On the boundary $\vert z \vert = 1$ of the unit circle, the series
(\ref{Ser_2f1}) diverges if $\mathrm{Re} (a+b-c) \ge 1$, it converges
\emph{absolutely} for $\mathrm{Re} (a+b-c) < 0$, and it converges
\emph{conditionally} for $0 \le \mathrm{Re} (a+b-c) < 1$, the point $z =
1$ being excluded.

Thus, a nonterminating hypergeometric series does not suffice for the
computation of the corresponding hypergeometric function ${}_2 F_1 (a,
b; c; z)$, which is in general a multivalued function defined in the
whole complex plane with branch points at $z = 1$ and $z =
\infty$. Instead, techniques which permit an analytic continuation from
the interior to the exterior of the unit circle are needed.

Sequence transformations can be used to accelerate the convergence of a
hypergeometric series ${}_2 F_1$ or to sum it in the case of divergence,
which corresponds to an analytic continuation. Let us for example
consider the following elementary special case of a hypergeometric
function ${}_2 F_1$ (Eq.\ (15.1.3) of \cite{AbSte72}):
\begin{eqnarray}
\ln (1+z) & = & \sum_{m=0}^{\infty} \, \frac {(-1)^m z^{m+1}} {m+1} 
\nonumber \\
& = & z \, {}_2 F_1 (1, 1; 2; -z) \, .
\label{HygT3_1}
\end{eqnarray}
The infinite series converges only for $\vert z \vert < 1$, whereas the
logarithm is with the exception of the cut along $- \infty < z \le - 1$
defined in the whole complex plane.

In Table \ref{Tab_3_1}, Wynn's epsilon algorithm, Eq.\ (\ref{eps_al}),
Brezinski's theta algorithm, Eq.\ (\ref{theta_al}), Levin's
transformation $d_{k}^{(n)} (\zeta, s_n)$, Eq.\ (\ref{dLevTr}), and the
closely related sequence transformation $\delta_{k}^{(n)} (\zeta, s_n)$,
Eq.\ (\ref{dSidTr}), are applied to the partial sums
\begin{equation}
s_n (z) \; = \; \sum_{m=0}^{n} \, \frac {(-1)^m z^{m+1}} {m+1}
\label{ParSumT3_1}
\end{equation}
of the hypergeometric series (\ref{HygT3_1}) for $z = 7/2$. The
approximations to the limit in Table \ref{Tab_3_1} were always chosen in
such a way that the transforms with the highest possible transformation
order were computed from a given set of input data. Thus, in the case of
the epsilon algorithm, the approximations to the limit were chosen
according to (\ref{EpsApprLim}), in the case of the theta algorithm,
they were chosen according to (\ref{ThetApproxLim}), and in the case of
the sequence transformations $d_{k}^{(n)} (\zeta, s_n)$ and
$\delta_{k}^{(n)} (\zeta, s_n)$, they were chosen according to
(\ref{LevTypeApprLim}).

The second column of Table \ref{Tab_3_1}, which displays the partial
sums (\ref{ParSumT3_1}), shows that the hypergeometric series
(\ref{HygT3_1}) for $\ln (1+z)$ diverges quite strongly for $z =
7/2$. Nevertheless, it is apparently possible to sum this divergent
series to its correct value.

The results in Table \ref{Tab_3_1} also show that Wynn's epsilon
algorithm, which in the case of a power series produces Pad\'e
approximants according to (\ref{Eps_Pade}), is contrary to a widespread
belief not necessarily the most powerful transformation. The theta
algorithm and in particular the two Levin-type transformations
$d_{k}^{(n)} (\zeta, s_n)$ and $\delta_{k}^{(n)} (\zeta, s_n)$, which
use the first term neglected in the partial sum as a remainder estimate
according to (\ref{dRemEst}), produce significantly better summation
results.

The other sequence transformations discussed in Appendix \ref{App_A}
give better results than Wynn's epsilon algorithm, but are less
effective than the transformations $d_{k}^{(n)} (\zeta, s_n)$ and
$\delta_{k}^{(n)} (\zeta, s_n)$. For example, the two approximations
with the highest possible transformation orders, which Aitken's iterated
$\Delta^2$ process (\ref{It_Aitken}) produces from the partial sums $s_0
(z)$, $s_1 (z)$, \ldots, $s_{15} (z)$ according to (\ref{AitApprLim}),
are
\begin{eqnarray}
{\cal A}_{7}^{(0)} & = & 1.504~077~397~173 \, ,  
\\
{\cal A}_{7}^{(1)} & = & 1.504~077~396~169 \, .
\end{eqnarray}  
Similarly, the iteration (\ref{thetit_1}) of Brezinski's theta algorithm
produces according to (\ref{ThetitApproxLim}) the approximants
\begin{eqnarray}
{\cal J}_{4}^{(2)} & = & 1.504~077~404~830 \, ,
\\
{\cal J}_{5}^{(0)} & = & 1.504~077~394~094 \, .
\end{eqnarray}

These results show that sequence transformations can be very
useful. Nevertheless, in the case of a \emph{real} real argument $z$ it
is actually not necessary to to use sequence transformations for doing
the analytic continuations or for speeding up convergence. Instead, it
is often simpler to exploit some known mathematical properties: Unless
certain linear combinations of the parameters $a$, $b$, and $c$ are
positive or negative integers, a hypergeometric function ${}_2 F_1 (a,
b; c; z)$ can be expressed as the sum of two other ${}_2 F_1$'s with a
transformed argument $w = 1 - z$, $w = 1/z$, $w = 1/(1-z)$, or $w =
1-1/z$, respectively. Thus, the argument $w$ of the two resulting
hypergeometric series can normally be chosen in such a way that the two
new series in $w$ either converge, if the original series in $z$
diverges, or that they converge more rapidly if the original series
converges too slowly to be numerically useful.

For example, if $\vert 1-z \vert < 1$ and if $c - a - b$ is not a
positive or negative integer, then we can use the analytic continuation
formula (Eq.\ (15.3.6) of \cite{AbSte72})
\begin{eqnarray}
\lefteqn{{}_2 F_1 (a, b; c; z)} \nonumber \\
& = &
\frac {\Gamma (c) \Gamma (c-a-b)} {\Gamma (c-a) \Gamma (c-b)} \,
{}_2 F_1 (a, b; a+b-c+1; 1-z) \nonumber \\
& & + \,
\frac {\Gamma (c) \Gamma (a+b-c)} {\Gamma (a) \Gamma (b)} \,
(1-z)^{c-a-b} \nonumber \\
& & \times \, {}_2 F_1 (c-a, c-b; c-a-b+1; 1-z) \, .
\label{AC_1}
\end{eqnarray}
Let us now assume that $z$ is only slightly smaller than 1. Then, the
convergence of the original hypergeometric series ${}_2 F_1 (a, b; c;
z)$ will be very bad. However, the two hypergeometric series on the
right-hand side with argument $1-z$ will converge rapidly in the
vicinity of $z = 1$.

With the help of this or similar analytic continuation formulas it is
normally possible to compute a hypergeometric function ${}_2 F_1 (a, b;
c; z)$ with \emph{real} argument $z$ effectively since it is possible to
find for every argument $z \in (- \infty, + \infty)$ two hypergeometric
series with an argument $\vert w \vert \le 1/2$ (see p.\ 127 of
\cite{Tem96} or Table I of \cite{For97}).

Unfortunately, this approach does not necessarily work in the case of
complex arguments $z$. Consider the points 
\begin{equation}
z_{1,2} \; = \; \frac{1 \pm \sqrt{3}}{2} \, , 
\end{equation}
which both lie on the boundary of the circle of convergence because of
$\vert z_{1,2} \vert = 1$. In practice, it is either impossible or not
feasible to evaluate a nonterminating hypergeometric series ${}_2 F_1$
by adding up its terms if its argument $z$ lies on the boundary of the
unit circle. Consequently, something has to be done to speed up
convergence or to accomplish a summation in the case of
divergence. Unfortunately, the analytic continuation formulas of the
type of (\ref{AC_1}) do not improve the situation if $z = z_{1,2}$ since
\begin{mathletters}
\begin{eqnarray}
1 - z_{1,2} & = & z_{2,1} \, , \\
1/z_{1,2} & = & z_{2,1} \, , \\
1/(1-z_{1,2}) & \; = \; & z_{1,2} \, , \\
1 - 1/z_{1,2} & \; = \; & z_{1,2} \, .
\end{eqnarray}
\end{mathletters}
Hence, the analytic continuation formulas cannot help if the $z$ is
close to $z_{1,2}$.

However, sequence transformations work for $z = z_{1,2}$. Let us
consider the following elementary special case of a hypergeometric
function ${}_2 F_1$ (p.\ 38 of \cite{MaObSo66}):
\begin{equation}
(1+z) (1-z)^{-2 \alpha-1} \; = \; 
{}_2 F_1 (2\alpha, \alpha+1; \alpha; z) \, .
\label{HygT3_2}
\end{equation}
For $\alpha > 0$, the hypergeometric series converges in the interior of
the unit circle, but it diverges on its boundary.

In Table \ref{Tab_3_2}, Wynn's epsilon algorithm, Eq.\ (\ref{eps_al}),
and Levin's transformation $d_{k}^{(n)} (\zeta, s_n)$, Eq.\
(\ref{dLevTr}), are used to sum the hypergeometric series
(\ref{HygT3_2}) with $\alpha = 1/3$ on the boundary of the unit circle,
i.e., they are applied to the partial sums
\begin{equation}
s_n (z) \; = \; 
\sum_{m=0}^{n} \, \frac{(2/3)_m (4/3)_m}{(1/3)_m m!} \, z^m 
\label{ParSumT3_2}
\end{equation}
with $z = z_1 = (1 + {\rm i} \sqrt{3})/2$.

The partial sums (\ref{ParSumT3_2}) in the second column of Table
\ref{Tab_3_2} display a very unusual sign pattern and grow slowly in 
magnitude. Nevertheless, both the epsilon algorithm as well as the Levin
transformation are apparently able to sum the hypergeometric series
(\ref{HygT3_2}) for $z = (1 + {\rm i} \sqrt{3})/2$.

The other sequence transformations discussed in Appendix \ref{App_A} are
apparently also able to sum the hypergeometric series (\ref{HygT3_2})
for $z = (1 + {\rm i} \sqrt{3})/2$. They give better results than Wynn's
epsilon algorithm, but are less effective than Levin's transformations
$d_{k}^{(n)} (\zeta, s_n)$. The two approximations with the highest
possible transformation orders, that can be produced from the partial
sums $s_0 (z)$, $s_1 (z)$, \ldots, $s_{15} (z)$ by Aitken's iterated
$\Delta^2$ process (\ref{It_Aitken}), by Brezinski's theta algorithm
(\ref{theta_al}) and its iteration (\ref{thetit_1}), and by the
Levin-type transformation (\ref{dSidTr}), are
\begin{eqnarray}
\lefteqn{{\cal A}_{7}^{(0)} \; = \;} \nonumber \\
& & \quad - \, 
1.113~340~798~057 \, + \, {\rm i} \, 1.326~827~896~288 \, , \\ 
\lefteqn{{\cal A}_{7}^{(1)} \; = \;} \nonumber \\
& & \quad - \, 
1.113~340~798~408 \, + \, {\rm i} \, 1.326~827~896~424 \, , \\ 
\lefteqn{\theta_{8}^{(2)} \; = \;} \nonumber \\
& & \quad - \, 
1.113~340~797~528 \, + \, {\rm i} \, 1.326~827~893~689 \, , \\ 
\lefteqn{\theta_{10}^{(0)} \; = \;} \nonumber \\
& & \quad - \, 
1.113~340~799~160 \, + \, {\rm i} \, 1.326~827~895~539 \, , \\ 
\lefteqn{{\cal J}_{4}^{(2)} \; = \;} \nonumber \\
& & \quad - \, 
1.113~340~798~249 \, + \, {\rm i} \, 1.326~827~894~967 \, , \\ 
\lefteqn{{\cal J}_{5}^{(0)} \; = \;} \nonumber \\
& & \quad - \, 
1.113~340~798~314 \, + \, {\rm i} \, 1.326~827~895~955 \, , \\ 
\lefteqn{\delta_{14}^{(0)} (1, s_0 (z)) \; = \;} \nonumber \\
& & \quad - \, 
1.113~340~798~314 \, + \, {\rm i} \, 1.326~827~895~955 \, , \\ 
\lefteqn{\delta_{15}^{(0)} (1, s_0 (z)) \; = \;} \nonumber \\
& & \quad - \, 
1.113~340~798~414 \, + \, {\rm i} \, 1.326~827~896~325 \, .
\end{eqnarray}

The numerical results presented here should suffice to support the claim
of the author that sequence transformations can be extremely useful
numerical tools for the evaluation of special functions in general
\cite{We89,We90,WeCi90,We94,We96d,JeMoSoWe99} and for the evaluation of
Gaussian hypergeometric series ${}_2 F_1$ in special. Moreover, the
numerical results shown above indicate that a computational algorithm,
which would be capable of evaluating of a hypergeometric function ${}_2
F_1$ with essentially arbitrary \emph{complex} argument $z$ and
parameters $a$, $b$, and $c$, should be a suitable combination of
analytic continuation formulas of the type of (\ref{AC_1}) with sequence
transformations. 

Here it should be taken into account that a hypergeometric series ${}_2
F_1$ depends on four essentially arbitrary complex quantities $a$, $b$,
$c$, and $z$. This makes the development of a {\em general\/} algorithm
for its computation difficult since a very large variety of special
cases and computationally different situations have to be taken into
account. Consequently, the development of such an algorithm would first
require extensive numerical studies which -- although undeniably
interesting -- would clearly be beyond the scope of this article.

\section{Hypergeometric Series with a Negative Third Parameter}
\label{sec:Ne3Pa2F1}

It is a direct consequence of its definition (\ref{Ser_2f1}) that a
hypergeometric series ${}_2 F_1 (a, b; c; z)$ terminates after a finite
number of terms if either $a$ or $b$ is a negative integer. Moreover,
${}_2 F_1 (a, b; c; z)$ makes sense if both $a$ and $c$ are negative
integers such that $a = - m$ and $c = -m - k$ with $k, m = 1, 2,
\ldots~$:
\begin{equation}
{}_2 F_1 (- m, b; -m - k; z) \; = \; \sum_{\mu=0}^{m} \,
\frac {(-m)_{\mu} (b)_{\mu}} {(-m-k)_{\mu} {\mu}!} \, z^{\mu} \, .
\end{equation}
If $c$ is a negative integer and if neither $a$ nor $b$ is a negative
integer, then it follows from (\ref{Ser_2f1}) that the hypergeometric
series ${}_2 F_1 (a, b; c; z)$ is in general undefined. However, the
following limit exists for $m = 0, 1, 2, \ldots$ (p.\ 38 of
\cite{MaObSo66}):
\begin{eqnarray}
\lefteqn{\lim_{c \to - m} \, \frac {1} {\Gamma (c)} \,
{}_2 F_1 (a, b; c; z) \; = \; 
\frac {(a)_{m+1} (b)_{m+1} z^{m+1}} {(m+1)!}} \nonumber \\
& & \quad \times \,  {}_2 F_1 (a+m+1, b+m+1; m+2; z) \, .
\end{eqnarray}

If $c$ is not a negative integer, there are no problems with Pochhammer
symbols in the denominators of the terms of the series that could become
zero. Nevertheless, unpleasant numerical problems occur even if $c$ is
just a nonintegral negative real number. These problems can be
demonstrated convincingly by trying to accelerate the convergence of two
hypergeometric series ${}_2 F_1 (a, b; c; z)$ which differ only by the
sign of the third parameter $c$.

For that purpose, we consider in Table \ref{Tab_4_1} and in Tables
\ref{Tab_4_2} and \ref{Tab_4_2_a}, respectively, the hypergeometric
series with $a = 3/7$, $b = 5/2$, $z = 77/100$, and $c = \pm 7/2$.

In Table \ref{Tab_4_1}, Wynn's epsilon algorithm, Eq.\ (\ref{eps_al}),
Brezinski's theta algorithm, Eq.\ (\ref{theta_al}), and the Levin-type
transformations $d_{k}^{(n)} (\zeta, s_n)$, Eq.\ (\ref{dLevTr}), and 
$\delta_{k}^{(n)} (\zeta, s_n)$, Eq.\ (\ref{dSidTr}), are applied to the 
partial sums
\begin{equation}
s_n (z) \; = \; 
\sum_{m=0}^{n} \, \frac{(3/7)_m (5/2)_m}{(7/2)_m m!} \, z^m
\label{ParSumT4_1}
\end{equation}
of the hypergeometric series ${}_2 F_1 (3/7, 5/2; 7/2; z)$ with $z =
77/100$. The results in Table \ref{Tab_4_1} show that the convergence of
this series can indeed be accelerated quite effectively by sequence
transformations. This is also true for the other sequence
transformations discussed in Appendix \ref{App_A}. The approximations
with the highest transformation orders, that can be obtained from the
partial sums $s_0 (z)$, $s_1 (z)$, \ldots, $s_{16} (z)$ by Aitken's
iterated $\Delta^2$ process (\ref{It_Aitken}) and by the iteration
(\ref{thetit_1}) of Brezinski's theta algorithm, are
\begin{eqnarray}
{\cal A}_{7}^{(1)} & = & 1.463~807~099~629 \, ,
\\
{\cal A}_{8}^{(0)} & = & 1.463~807~099~563 \, .
\\
{\cal J}_{5}^{(0)} & = & 1.463~807~143~254 \, ,
\\
{\cal J}_{5}^{(1)} & = & 1.463~807~103~421 \, .
\end{eqnarray}

In Tables \ref{Tab_4_2} and \ref{Tab_4_2_a}, we now consider a
hypergeometric series ${}_2 F_1$ which is identical with the one in
Table \ref{Tab_4_1} except that its third parameter is negative, i.e.,
we now have $c = -7/2$. Thus, in Table \ref{Tab_4_2} we use Aitken's
iterated $\Delta^2$ process, Eq.\ (\ref{It_Aitken}), Wynn's epsilon
algorithm, Eq.\ (\ref{eps_al}), and Brezinski's theta algorithm, Eq.\
(\ref{theta_al}), for the
acceleration of the convergence of the partial sums
\begin{equation}
s_n (z) \; = \; 
\sum_{m=0}^{n} \, \frac{(3/7)_m (5/2)_m}{(-7/2)_m m!} \, z^m \, .
\label{ParSumT4_2}
\end{equation}
of the hypergeometric series ${}_2 F_1 (3/7, 5/2; - 7/2; z)$ with $z =
77/100$, and in Table \ref{Tab_4_2_a} we use the iteration of
Brezinski's theta algorithm, Eq.\ (\ref{thetit_1}), and the Levin-type
transformations $d_{k}^{(n)} (\zeta, s_n)$, Eq.\ (\ref{dLevTr}), and
$\delta_{k}^{(n)} (\zeta, s_n)$, Eq.\ (\ref{dSidTr}).

So far, the transformation results had always been very good as well as
very reliable. In contrast, the results in Tables \ref{Tab_4_2} and
\ref{Tab_4_2_a} are very bad and not reliable at all. For small
transformation orders $n$, all transformations produce results which are
by 5 orders of magnitude too small, and in the case of Brezinski's theta
algorithm, the wrong results even seem to have converged with an
accuracy of 4 decimal digits. For increasing transformation orders $n$,
there occur sudden and unmotivated sign changes, and only if $n$
approaches 30, at least the epsilon algorithm in Table \ref{Tab_4_2} and
the Levin transformation $d_{n}^{(0)}$ in Table \ref{Tab_4_2_a} converge
to the correct result. For $n = 30$, the other transformations show no
indication of convergence and produce results which are still by some
orders of magnitude too small.

How can the disturbingly bad performance of sequence transformations in
Tables \ref{Tab_4_2} and \ref{Tab_4_2_a} be explained. The terms of a
hypergeometric series ${}_2 F_1 (a, b; c; z)$ satisfy the 2-term
recursion
\begin{eqnarray}
\lefteqn{\frac{(a)_{n+1} (b)_{n+1} z^{n+1}}{(c)_{n+1} (n+1)!}} 
\nonumber \\
&& \qquad \; = \; 
\frac{(a+n)(b+n)z}{(c+n)(n+1)} \, \frac{(a)_n (b)_n z^n}{(c)_n n!} \, .
\end{eqnarray}
Obviously, the factor $(a+n)(b+n)z/[(c+n)(n+1)]$ on the right-hand side
determines whether the terms increase or decrease in magnitude with
increasing $n$. As long as this factor is greater than one in magnitude,
the terms increase with increasing $n$, and as soon as this factor is
smaller than one, the terms decrease.
  
Thus, we only have to determine those values of $n$ which satisfy
\begin{equation}
\left\vert \frac{(a+n) (b+n) z}{(c+n) (n+1)} \right\vert \; = \; 1
\label{CfEq1}
\end{equation}
for given $a$, $b$, $c$, and $z$ in order to find out for which values
of the index $n$ the terms of a hypergeometric series ${}_2 F_1 (a, b;
c; z)$ change their growth pattern.

In the case of the hypergeometric series in Table \ref{Tab_4_1} with $a
= 3/7$, $b = 5/2$, $c = 7/2$, and $z = 77/100$, there is no $n > 0$
which satisfies condition (\ref{CfEq1}). This implies that the terms of
this series decrease monotonously in magnitude with increasing $n \ge
0$. Moreover, the results in Table \ref{Tab_4_1} show that the
convergence of this series can be accelerated quite effectively.
 
In the case of the hypergeometric series in Tables \ref{Tab_4_2} and
\ref{Tab_4_2_a} with $c = - 7/2$, the situation is more complicated
since condition (\ref{CfEq1}) is satisfied by $n \approx 22$. Thus, the
terms of this series initially increase up to $n = 22$, and only for $n
> 22$ they decrease and ultimately produce a convergent
result. Accordingly, the partial sums (\ref{ParSumT4_2}) of this
hypergeometric series initially look like the partial sums of a mildly
divergent series, and only for $n > 22$, they behave like the partial
sums of a convergent series. Therefore, it should not be too surprising
that sequence transformations perform poorly if they use as input data
only the partial sums (\ref{ParSumT4_2}) with $n \le 22$. However, even
for $22 \le n \le 30$, only the epsilon algorithm and the Levin
transformation $d_{n}^{(0)}$ ultimately converge to the correct
result. This provides strong evidence that the irregular input data with
small indices $n$ have a detrimental effect on the performance of
sequence transformations with large transformation orders, which also
use input data with a correct behavior.

We can test the hypothesis, that the irregular behavior of the initial
partial sums (\ref{ParSumT4_2}) with $n \le 22$ leads to the poor
performance of sequence transformations in Tables \ref{Tab_4_2} and
\ref{Tab_4_2_a}, by skipping the terms up to $n = 22$ in the
transformation processes. Accordingly, Wynn's epsilon algorithm, Eq.\
(\ref{eps_al}), and the Levin-type transformations $d_{k}^{(n)} (\zeta,
s_n)$, Eq.\ (\ref{dLevTr}), and $\delta_{k}^{(n)} (\zeta, s_n)$, Eq.\
(\ref{dSidTr}), use in in Table \ref{Tab_4_3} the modified partial sums
\begin{equation}
s_{n}^{(22)} (z) \; = \; 
\sum_{m=0}^{n+22} \, \frac{(3/7)_m (5/2)_m}{(-7/2)_m m!} \, z^m 
\label{ParSumT4_3}
\end{equation}
of the hypergeometric series ${}_2 F_1 (3/7, 5/2; - 7/2; z)$ with $z =
77/100$ as input data. This approach is possible since all sequence
transformations considered in this article are quasi-linear, i.e., they
satisfy (\ref{quasi_lin}).

The results in Table \ref{Tab_4_3} indeed confirm our hypothesis since
they are nearly as good as the results in Table \ref{Tab_4_1}, at least
with respect to the transformation orders that are needed to achieve a
given relative accuracy. 

The other sequence transformations discussed in Appendix \ref{App_A}
produce results that are less good than those shown in Table
\ref{Tab_4_3}. The two approximations with the highest possible
transformation orders, that can be produced from the partial sums
$s_0^{(22)} (z)$, $s_1^{(22)} (z)$, \ldots, $s_{20}^{(22)} (z)$ by
Aitken's iterated $\Delta^2$ process (\ref{It_Aitken}), by Brezinski's
theta algorithm (\ref{theta_al}) and its iteration (\ref{thetit_1}), are
\begin{eqnarray}
{\cal A}_{9}^{(1)} & = & 1.010~147~722~439 \cdot 10^{+5} \, ,
\\
{\cal A}_{10}^{(0)} & = & 1.010~147~537~701 \cdot 10^{+5} \, .
\\
\theta_{12}^{(1)} & = & 1.011~462~051~628 \cdot 10^{+5} \, ,
\\
\theta_{12}^{(2)} & = & 1.011~462~011~501 \cdot 10^{+5} \, ,
\\
{\cal J}_{7}^{(0)} & = & 1.010~233~908~825 \cdot 10^{+5} \, ,
\\
{\cal J}_{7}^{(1)} & = & 1.010~176~054~786 \cdot 10^{+5} \, .
\end{eqnarray}

The results in Table \ref{Tab_4_3} show that a hypergeometric series
${}_2 F_1$ with a negative third parameter can be evaluated reliably
with the help of sequence transformations if the nonregular leading
terms are excluded from the transformation processes. Unfortunately,
this may lead to new problems since the number of terms, that initially
grow in magnitude, may become quite large, in particular if $z$ is close
to one. For example, if we increase the argument of the hypergeometric
series in Table \ref{Tab_4_3} from $z = 77/100$ to $z = 87/100$ or to $z
= 97/100$, then the number of terms, which initially grow in magnitude
and have to be skipped, grow from $n = 22$ to $n = 40$ or even to $n =
179$. Moreover, more negative values of $c$ also increase the number of
terms that have to be skipped. For instance, if we consider the
hypergeometric series ${}_2 F_1 (3/7, 5/2; - 13/2; z)$ with $z =
77/100$, $z = 87/100$, or $z = 97/100$, then we obtain $n = 35$, $n =
63$, or $n = 279$, respectively. These examples show that it is in
principle possible to construct hypergeometric series ${}_2 F_1$ which
can only be evaluated reliably with the help of sequence transformations
if a very large number of terms is skipped in the transformation
process.

In the case of a Gaussian hypergeometric series, this poses no
unsurmountable problems. Firstly, it is a triviality to compute the
terms, even for very large indices. Consequently, a brute force
evaluation of a Gaussian hypergeometric series by adding up the terms is
possible as long as the argument is not too close to the boundary of the
circle of convergence. Secondly, the highly developed mathematical
theory of these functions makes it possible to simplify the numerical
task with the help of known transformation formulas. For example, if the
number of terms, that have to be skipped, is large because the argument
$z$ of the hypergeometric series with a negative third parameter is
close to one, then it may be a good idea to use the analytic
continuation formula (\ref{AC_1}). This would not necessarily solve the
principal problems due to a negative third parameter, but the argument
$1-z$ of the two new hypergeometric series would then be small.

However, the probably simplest approach would be the use of recurrence
formulas. In this approach, one would have to evaluate two
hypergeometric series ${}_2 F_1$ with suitable positive values of the
third parameter, for example with the help of sequence transformations,
and to compute recursively the numerical value of desired hypergeometric
series with a negative third parameter. Of course, many recurrence
formulas are known. Nevertheless, numerous new three-term recurrence
formulas satisfied by the Gaussian hypergeometric function ${}_2 F_1 (a,
b; c; z)$ are derived in Appendix \ref{App_C}.

Unfortunately, these alternative approaches are in general not available
if we have to evaluate an infinite series whose terms are determined
numerically and only behave like the terms of a hypergeometric series
${}_2 F_1 (a, b; c; z)$ with a negative third parameter. In such a case,
it cannot be excluded that the number of terms, that have to be skipped,
would be so large that their computation would no longer be
feasible. Then, neither the conventional process of successively adding
up the terms of the series nor sequence transformations would be able to
provide reliable approximations to the value of such an infinite
series. 

\section{Summary and Conclusions}
\label{sec:SumConclu}

A sequence transformation is a rule ${\cal T}$ which transforms a
sequence $\{ s_n \}_{n=0}^{\infty}$ with the (generalized) limit $s$ to
another sequence $\{ s'_n \}_{n=0}^{\infty}$ having the same
(generalized) limit $s$ but different remainders $r'_n = s'_n - s$. The
transformation process was successful if the new sequence has better
numerical properties than the original sequence. For example, a sequence
transformation ${\cal T}$ accelerates convergence, if the transformed
remainders $r'_n = s'_n -s$ converge more rapidly than the original
remainders $r_n = s_n - s$ according to (\ref{CondConvAcc}), and ${\cal
T}$ sums a divergent sequence, whose remainders $r_n$ do not vanish as
$n \to \infty$, to its generalized limit $s$ if the transformed
remainders $r'_n$ approach zero as $n \to \infty$.

As discussed in Section \ref{sec:Paths}, a sequence transformation
${\cal T}$ can be represented by an infinite set of doubly indexed
quantities $T_{k}^{(n)}$ with $k, n \ge 0$ that can be displayed in a
two-dimensional array called the table of ${\cal T}$. The superscript
$n$ denotes the minimal index occurring in the finite string $\{ s_n,
s_{n+1}, \ldots s_{n+l} \}$ of sequence elements used for the
computation of a given $T_k^{(n)}$. The subscript $k$ -- usually called
the order of the transformation -- is a measure for the complexity of
the transformation process which yields $T_k^{(n)}$.

A convergence acceleration or summation process tries to obtain a better
approximation to the (generalized) limit of the input sequence by
proceeding on a certain path in the table of ${\cal T}$. There is an in
principle unlimited variety of different paths, but in practical
applications either order-constant paths defined in (\ref{OrdConstPath})
or index-constant paths defined in (\ref{IndConstPath}) are normally
used.

On an order-constant path, a set $\{ s_n, s_{n+1}, \ldots s_{n+l} \}$ of
input data of fixed length is used and the starting index $n$ of this
set is increased successively. In contrast, an index-constant path uses
sets of input data of increasing length, and it is always tried to
compute from a given set of input data that element $T_k^{(n)}$ which
has the highest possible transformation order $k$.

Order-constant and index-constant paths differ substantially. For
example, on an index-constant path the available information is exploited
more efficiently than on an order-constant path. Accordingly,
index-constant paths normally produce better transformation
results. Moreover, order-constant paths cannot be used for the summation
of a divergent sequence since increasing $n$ in the set $\{ s_n,
s_{n+1}, \ldots s_{n+l} \}$ of input data normally only increases
divergence.

If, however, the leading elements of the input sequence $\{ s_n
\}_{n=0}^{\infty}$ behave irregularly, the principal advantages of 
index-constant paths can easily turn into disadvantages: If the sets $\{
s_n, s_{n+1}, \ldots s_{n+l} \}$ of input data with increasing $l$ have
a sufficiently small starting index $n$, then all transforms
$T_{k}^{(n)}$ will be affected by irregular input data, albeit to a
different degree. As shown by Tables \ref{Tab_4_2} and \ref{Tab_4_2_a},
this can lead to unreliable or even completely nonsensical
transformation results. In such a case, it is necessary to exclude the
irregular input data from the transformation process. This can be
accomplished by using either an order-constant path or an index-constant
path with a sufficiently large starting index.

In Section \ref{sec:2F1}, the Gaussian hypergeometric function ${}_2 F1
(a, b; c; z)$ is discussed, which is in general a multivalued function
defined in the whole complex plane with branch points at $z = 1$ and
$\infty$. However, it is defined by the power series (\ref{Ser_2f1})
which only converges for $\vert z \vert < 1$. Accordingly, sequence
transformations can either be used for speeding up convergence or for
accomplishing an analytic continuation in the case of divergence. In
Table \ref{Tab_3_1}, an alternating hypergeometric series for $\ln
(1+z)$ is summed effectively by sequence transformations for an
argument $z = 7/2$ which is is far away from the unit circle, and in
Table \ref{Tab_3_2}, another special hypergeometric series, which
converges only in the interior of the unit circle, is evaluated for an
argument $z = (1+{\rm i}\sqrt{3})/2$ that is located on the boundary of
the unit chicle.

For most parameters $a$, $b$, and $c$, sequence transformations greatly
facilitate the evaluation of a hypergeometric series ${}_2 F_1$, and the
good transformation results presented in Tables \ref{Tab_3_1} and
\ref{Tab_3_2} are fairly typical. However, as discussed in Section
\ref{sec:Ne3Pa2F1}, there is an important and instructive exception: If
the third parameter $c$ of a hypergeometric series is a negative real
number, then the terms of this series initially increase in magnitude
and look even for $\vert z \vert < 1$ like the terms of a mildly
divergent series. Only for sufficiently large values of the index, the
terms decrease and ultimately produce a convergent result. 

The use of these irregular terms as input data seriously affects the
performance of sequence transformations and leads to unreliable and
sometimes even completely nonsensical results.  This is demonstrated by
applying in Table \ref{Tab_4_1} and in Tables \ref{Tab_4_2} and
\ref{Tab_4_2_a}, respectively, sequence transformations to the partial
sums of the hypergeometric series ${}_2 F_1 (a, b; c; z)$ with $a =
3/7$, $b = 5/2$, $z = 77/100$, and $c = \pm 7/2$.

In Table \ref{Tab_4_1}, sequence transformations are applied to the
hypergeometric series with the positive value of the third parameter. As
expected, the transformation results are very good. However, in Tables
\ref{Tab_4_2} and \ref{Tab_4_2_a}, where the hypergeometric series with
the negative value of the third parameter is considered, the
transformation results are both unreliable and bad.

Nevertheless, it is possible to compute the hypergeometric series ${}_2
F_1 (3/7, 5/2; - 7/2; 77/100)$ efficiently and reliably with the help of
sequence transformations. However, one cannot use an index-constant path
with a small minimal index, as it was done in Tables \ref{Tab_4_2} and
\ref{Tab_4_2_a}. Instead, one should either use an order-constant path
or an index-constant path with a sufficiently large minimal index, as it
was done in Table \ref{Tab_4_3}, where all irregular terms were excluded
from the transformation processes.

The sequence transformations, which are used in this article, are all
described in Appendix \ref{App_A}. The use of several transformations
was quite intentional. The author wanted to make clear that problems due
to irregular input data are not restricted to some special sequence
transformations only. Of course, the results in Tables \ref{Tab_4_2} and
\ref{Tab_4_2_a} show that different sequence transformations respond 
differently to irregular input data. However, this is quite helpful and
can protect us against misinterpretations. For example, in Table
\ref{Tab_4_2} Brezinski's theta algorithm produced transformation
results which seemed to have converged with an accuracy of 4 decimal
digits, but were actually by 5 orders of magnitude too
small. Fortunately, the other transformations in Tables \ref{Tab_4_2}
and \ref{Tab_4_2_a} produced different results, which provided strong
evidence that the transformation results were unreliable. Consequently,
it is recommendable to use in convergence acceleration and summation
processes more than a single transformation whenever possible. This is
particularly important if numerically determined data are to be
transformed, about which very little is known. If several different
sequence transformations produce consistent results, then it is very
likely that these results are indeed correct although it is of course
clear that purely numerical results cannot be a substitute for a rigorous
mathematical proof.

It looks like a contradiction that in Appendix \ref{App_B} only the
Levin-type transformation (\ref{dSidTr}) was used for the computation of
the infinite coupling limit $k_3$ of the sextic anharmonic
oscillator. However, the divergent perturbation series (\ref{rPT_k3}),
whose leading coefficients also show an irregular behavior, constitutes
a very demanding summation problem, for which the other transformations
discussed in Appendix \ref{App_A} are not powerful enough.

Here, it should be emphasized once more that the two examples considered
in this article -- the Gaussian hypergeometric series ${}_2 F_1 (a, b;
c; z)$ with a negative third parameter and the divergent perturbation
series (\ref{rPT_k3}) for the infinite coupling limit $k_3$ -- are
comparatively simple model problems, and the leading irregular terms of
their series expansions pose no unsurmountable computational
problems. This is largely due to the fact that it is relatively easy to
find out which terms behave irregularly. In the case of the perturbation
series (\ref{rPT_k3}), we are in the fortunate situation that the
leading large-$n$ asymptotics (\ref{cAs3}) of the coefficients
$c_{n}^{(3)}$ is known, and in the case of the hypergeometric series,
one only has to solve (\ref{CfEq1}) in order to find out for which
indices $n$ the terms change their growth pattern. Moreover, in the case
of a ${}_2 F_1$ there are numerous alternative computational
approaches. For example, with the help of recurrence formulas the
evaluation of Gaussian hypergeometric series with a negative third
parameter can be avoided completely. Consequently, in Appendix
\ref{App_C} numerous new recurrence formulas satisfied by a ${}_2 F_1
(a, b; c; z)$ are derived.

Finally, the author wishes to express his hope that this article will
inspire additional research on the evaluation of special functions with
the help of sequence transformations. There can be no doubt that
sequence transformations are normally extremely useful tools for the
evaluation of special functions, and since the terms of the series
expansions for special functions are comparatively simple and explicitly
known, we can even hope to gain additional insight from those cases in
which sequence transformations fail to produce good transformation
results.

\acknowledgments

The author thanks the Fonds der Chemischen Industrie for
financial support.

\appendix

\section{Sequence Transformations}
\label{App_A}

This appendix gives a short description of all the sequence
transformations that are used in this article. Further details plus
additional references can be found in \cite{BreRZa91,Wi81,We89}. Here,
the same notation as in \cite{We89} is used.

One of the oldest sequence transformations (see for instance pp.\ 90 -
91 of \cite{Br91}) is Aitken's $\Delta^2$ formula \cite{Ait26}:
\begin{equation}
{\cal A}_1^{(n)} \; = \; s_n \; - \; \frac
{[ \Delta s_n ]^2} {\Delta^2 s_n} \, .
\label{Aitken}
\end{equation}
The (forward) difference operator $\Delta$ acts for all integers $n \ge
0$ on a function $f (n)$ according to
\begin{equation}
\Delta f (n) \; = \; f (n+1) - f (n) \, .
\end{equation}
The Aitken formula (\ref{Aitken}) is by construction exact for model
sequences of the type $s_n = s + c \lambda^n$ with $c \ne 0$ and
$\lambda \ne 1$. If the numerical values of three consecutive elements
$s_{n}$, $s_{n+1}$, and $s_{n+2}$ of this model sequence are known, then
the (generalized) limit $s$ of this sequence can be computed according
to 
\begin{equation}
{\cal A}_1^{(n)} \; = \; s \, ,
\end{equation}
no matter whether the sequence converges ($\vert \lambda \vert < 1$) or
diverges ($\vert \lambda \vert > 1$).

The power and practical usefulness of Aitken's $\Delta^2$ formula is of
course limited since it is designed to eliminate only a single
exponential term from the elements of the model sequence mentioned
above. However, the quantities ${\cal A}_1^{(n)}$ can again be used as
input data in (\ref{Aitken}). Hence, the $\Delta^2$ process can be
iterated, yielding the following nonlinear recursive scheme
\cite[Eq.\ (5.1-15)]{We89}:
\begin{mathletters}
\label{It_Aitken} 
\begin{eqnarray}
\lefteqn{} \nonumber \\
{\cal A}_0^{(n)} & = & s_n \, , \\
{\cal A}_{k+1}^{(n)} & = & {\cal A}_{k}^{(n)} -
\frac
{\bigl[\Delta {\cal A}_{k}^{(n)}\bigr]^2}
{\Delta^2 {\cal A}_{k}^{(n)}} \, .
\end{eqnarray}
\end{mathletters}%
In this article, the difference operator $\Delta$ acts only on the
superscript $n$ and not on the subscript $k$ of a doubly indexed
quantity like ${\cal A}_{k}^{(n)}$. A more detailed discussion of
Aitken's iterated $\Delta^2$ process as well as additional references
can for instance be found in Section 5 of \cite{We89} or in
\cite{We2000}. The iteration of other sequence transformations is
discussed in \cite{We91}.

In the case of Aitken's iterated $\Delta^2$ process, the approximation
to the limit of the input sequence with the highest possible
transformation order depends upon the index $m$ of the last sequence
element $s_m$ which was used in the recursion. If $m$ is either even or
odd, $m = 2 \mu$ or $m = 2 \mu + 1$, respectively, the approximations to
the limit are chosen according to
\begin{eqnarray}
\{ s_0, s_1, \ldots , s_{2 \mu}\} & \longrightarrow &
{\cal A}_{\mu}^{(0)} \, ,
\\
\{ s_1, s_2, \ldots , s_{2 \mu + 1}\} & \longrightarrow &
{\cal A}_{\mu}^{(1)} \, .
\end{eqnarray}
As in the case of Wynn's $\epsilon$ algorithm, these two relationships
can with the help of the notation $\Ent {x}$ for the integral part of
$x$ be expressed by a single equation (Eq.\ (5.2-6) of \cite{We89}):
\begin{eqnarray}
\lefteqn{\left\{ s_{m - 2 \Ent {m/2}}, s_{m - 2 \Ent {m/2} + 1},
\ldots , s_m \right\} \qquad \quad} \nonumber \\
& \qquad \quad \longrightarrow &
{\cal A}_{\Ent {m/2}}^{(m - 2 \Ent {m/2})} \, .
\label{AitApprLim}
\end{eqnarray}

The behavior of many practically relevant convergent sequences $\{ s_{n}
\}_{n=0}^{\infty}$ can be characterized by the asymptotic condition
\begin{equation}
\lim_{n \to \infty} \, \frac{s_{n+1} - s}{s_{n} - s} \; = \; \rho \, ,
\label{AsyConvCond}
\end{equation}
where $s = s_{\infty}$ is the limit of the sequence $\{ s_n
\}_{n=0}^{\infty}$. This condition closely resembles the well known
ratio test for infinite series. A convergent sequence satisfying
(\ref{AsyConvCond}) with $\vert \rho \vert < 1$ is called {\em linearly}
convergent, and it is called {\em logarithmically} convergent if $\rho =
1$.

It is one of the major weaknesses of the otherwise very powerful and
very useful epsilon algorithm (\ref{eps_al}) that it does not work in
the case of logarithmic convergence. Brezinski showed that this
principal weakness can be overcome by a suitable modification of the
recursive scheme (\ref{eps_al}), which leads to the so-called theta
algorithm
\cite{Bre71}: 
\begin{mathletters}
\label{theta_al}
\begin{eqnarray}
\theta_{-1}^{(n)} & = & 0 \, , \qquad
\theta_0^{(n)} \; = \; s_n \, , \\
\theta_{2 k + 1}^{(n)} & = & \theta_{2 k-1}^{(n+1)}
\, + \, 1 / [\Delta \theta_{2 k}^{(n)}] \, , \\
\theta_{2 k+2}^{(n)} & = & \theta_{2 k}^{(n+1)} \, +
\, \frac
{[\Delta \theta_{2 k}^{(n+1)}] \,
[\Delta \theta_{2 k + 1}^{(n+1)}]}
{\Delta^2 \theta_{2 k+1}^{(n)}} \, .
\end{eqnarray}
\end{mathletters}%
\noindent
As in the case of Aitken's iterated $\Delta^2$ process
(\ref{It_Aitken}), it is assumed that the difference operator $\Delta$
acts only upon the superscript $n$ and not on the subscript $k$.

Again, the approximation to the limit of the input sequence depends
upon the index $m$ of the last sequence element $s_m$ which was used in
the recursion. If we have $m = 3 \mu$, $m = 3 \mu + 1$, or $m = 3 \mu +
2$, respectively, the approximations to the limit with the highest
transformation orders are chosen according to
\begin{eqnarray}
\{ s_0, s_1, \ldots ,s_{3 \mu}\} & \longrightarrow &
\theta_{2 \mu}^{(0)} \, ,
\\
\{ s_1, s_2, \ldots ,s_{3 \mu + 1}\} & \longrightarrow &
\theta_{2 \mu}^{(1)} \, ,
\\
\{ s_2, s_3, \ldots ,s_{3 \mu + 2}\} & \longrightarrow &
\theta_{2 \mu}^{(2)} \, .
\end{eqnarray}
These three relationships can be expressed by a single equation
(Eq.\ (10.2-8) of \cite{We89}):
\begin{eqnarray}
\lefteqn{\left\{ s_{m - 3 \Ent {m/3}}, s_{m - 3 \Ent {m/3} + 1}, \ldots ,
s_m \right\} \qquad \quad} \nonumber \\
& \qquad \quad \longrightarrow &
\theta_{2 \Ent {m/3}}^{(m - 3 \Ent {m/3})} \, .
\label{ThetApproxLim}
\end{eqnarray}
Further details on the theta algorithm as well as additional references
can be found in Section 2.9 of \cite{BreRZa91} or in Sections 10
and 11 of \cite{We89}.

As for example discussed in \cite{We91}, new sequence transformations
can be constructed by iterating explicit expressions for sequence
transformations with low transformation orders. The best known example
of such an iterated sequence transformation is probably Aitken's
iterated $\Delta^2$ process (\ref{It_Aitken}) which is obtained by
iterating (\ref{Aitken}).

The same approach is also possible in the case of the theta algorithm
(\ref{theta_al}). A suitable closed-form expression, which may be
iterated, is (Eq.\ (10.3-1) of \cite{We89})
\begin{eqnarray}
\lefteqn{\vartheta_2^{(n)} \; = \; s_{n+1}} \nonumber \\
&& \quad - \, \frac
{\bigl[\Delta s_n\bigr] \bigl[\Delta s_{n+1}\bigr] 
\bigl[\Delta^2 s_{n+1}\bigr]}
{\bigl[\Delta s_{n+2}\bigr] \bigl[\Delta^2 s_n\bigr] -
\bigl[\Delta s_n\bigr] \bigl[\Delta^2 s_{n+1}\bigr]} \, .
\label{theta_2_1}
\end{eqnarray}
The iteration of this expression yields the following nonlinear
recursive scheme (Eq.\ (10.3-6) of \cite{We89}):
\begin{mathletters}
\label{thetit_1}
\begin{eqnarray}
&& {\cal J}_0^{(n)} \; = \; s_n \, ,
\\
&& {\cal J}_{k+1}^{(n)} \; = \;
{\cal J}_k^{(n+1)} \, - \, \nonumber \\ 
&& \quad \frac
{\bigl[ \Delta {\cal J}_k^{(n)} \bigr] 
\bigl[ \Delta {\cal J}_k^{(n+1)} \bigr]
\bigl[ \Delta^2 {\cal J}_k^{(n+1)} \bigr]}
{\bigl[ \Delta {\cal J}_k^{(n+2)} \bigr] 
\bigl[ \Delta^2 {\cal J}_k^{(n)} \bigr]  - 
\bigl[ \Delta {\cal J}_k^{(n)} \bigr] 
\bigl[\Delta^2 {\cal J}_k^{(n+1)} \bigr]} \, .
\end{eqnarray}
\end{mathletters}
In convergence acceleration and summation processes, the iterated
transformation ${\cal J}_k^{(n)}$ has similar properties as Brezinski's
theta algorithm from which it was derived: Both transformations are very
powerful as well as very versatile. ${\cal J}_k^{(n)}$ is not only an
effective accelerator for linear convergence as well as able to sum
divergent alternating series, but it is also able to accelerate the
convergence of many logarithmically convergent sequences and series
\cite{We89,We91,BhoBhaRoy89,Sab87,Sab91,Sab92,Sab95}.

Again, the approximation to the limit of the input sequence depends
upon the index $m$ of the last sequence element $s_m$ which was used in
the recursion. If we have $m = 3 \mu$, $m = 3 \mu + 1$, or $m = 3 \mu +
2$, respectively, the approximations to the limit with the highest
transformation orders are chosen according to
\begin{eqnarray}
\{ s_0, s_1, \ldots ,s_{3 \mu}\} & \longrightarrow &
{\cal J}_{\mu}^{(0)} \, ,
\\
\{ s_1, s_2, \ldots ,s_{3 \mu + 1}\} & \longrightarrow &
{\cal J}_{\mu}^{(1)} \, ,
\\
\{ s_2, s_3, \ldots ,s_{3 \mu + 2}\} & \longrightarrow &
{\cal J}_{\mu}^{(2)} \, .
\end{eqnarray}
These three relationships can be expressed by a single equation
(Eq.\ (10.4-7) of \cite{We89}):
\begin{eqnarray}
\lefteqn{\left\{ s_{m - 3 \Ent {m/3}}, s_{m - 3 \Ent {m/3} + 1}, 
\ldots , s_m \right\}} \nonumber \\
& \qquad \quad \longrightarrow &
{\cal J}_{\Ent {m/3}}^{(m - 3 \Ent {m/3})} \, .
\label{ThetitApproxLim}
\end{eqnarray}

So far, only sequence transformations were considered which use as input
data the elements of the sequence to be transformed. However, in some
cases structural information on the dependence of the remainders $r_n$
on the index $n$ is available. For example, it is well known that the
truncation error of a convergent series with strictly alternating and
monotonously decreasing terms is bounded in magnitude by the first term
not included in the partial sum and has the same sign as this term (see
for instance p.\ 132 of \cite{Kn64}). The first term neglected is also
the best simple estimate for the truncation error of a strictly
alternating nonterminating hypergeometric series ${}_2 F_0 (\alpha,
\beta; - x)$ with $\alpha, \beta, x > 0$ (Theorem 5.12-5 of 
\cite{Ca77}). Such an information on the behavior of the truncation 
errors can be extremely helpful in a convergence acceleration or
summation process. Unfortunately, the sequence transformations
considered so far are not able to benefit from it.

A convenient way of incorporating such an information into the
transformation process consists in the use of remainder estimates $\{
\omega_n \}_{n=0}^{\infty}$. Because of the explicit incorporation of the
information contained in the remainder estimates, sequence
transformations of that kind are potentially very powerful and as well
as very versatile.

The best-known example of such a sequence transformation is Levin's
transformation \cite{Le73} which is both very versatile and very
powerful \cite{BreRZa91,We89,SmFo79,SmFo82,Hom2000}:
\begin{eqnarray}
\lefteqn{{\cal L}_{k}^{(n)} (\zeta, s_n, \omega_n)}
\nonumber \\
& & \; = \; \frac
{\displaystyle
\sum_{j=0}^{k} \; ( - 1)^{j} \; {{k} \choose {j}} \;
\frac
{(\zeta + n +j )^{k-1}} {(\zeta + n + k )^{k-1}} \;
\frac {s_{n+j}} {\omega_{n+j}} }
{\displaystyle
\sum_{j=0}^{k} \; ( - 1)^{j} \; {{k} \choose {j}} \;
\frac
{(\zeta + n +j )^{k-1}} {(\zeta + n + k )^{k-1}} \;
\frac {1} {\omega_{n+j}} } \; .
\label{LevTr}
\end{eqnarray}
The shift parameter $\zeta$ has to be positive in order to admit $n = 0$
in (\ref{LevTr}). The most obvious choice, which is always used in this
article, is $\zeta = 1$. Recurrence formulas for the numerator and
denominator sums of ${\cal L}_{k}^{(n)} (\zeta, s_n, \omega_n)$ can be
found in Section 7.2 of \cite{We89}.

Levin's transformation is based on the assumption that the remainders
$r_n$ of the input sequence can for all $n \ge 0$ be approximated by a
remainder estimate $\omega_n$, which should be chosen such that $s_n - s
= \omega_n \bigl[ c + O(n^{-1}) \bigr]$ as $n \to \infty$, multiplied by
a polynomial in $1/(n+\zeta)$ with $\zeta > 0$. Levin \cite{Le73}
introduced several simple remainder estimates for infinite series which
give rise to several variants of Levin's sequence
transformation. Further details on Levin's transformation can for
instance be found in Section 7 of \cite{We89}.

In this article, we only consider the remainder estimate
\begin{equation}
\omega_n \; = \; \Delta s_n \; = \; a_{n+1} \, ,
\label{dRemEst}
\end{equation}
which was first proposed by Smith and Ford \cite{SmFo79}. It yields the
following variant of Levin's transformation (Eq.\ (7.3-9) of
\cite{We89}):
\begin{equation}
d_{k}^{(n)} (\zeta, s_n) \; = \; 
{\cal L}_{k}^{(n)} (\zeta, s_n, \Delta s_n) \, .
\label{dLevTr}
\end{equation}

Levin's transformation is based on the implicit assumption that the
ratio $[s_n - s]/\omega_n$ can be expressed as a power series in
$1/(n+\zeta)$. A different class of sequence transformations can be
derived by assuming that the ratio $[s_n - s]/\omega_n$ can be expressed
as a so-called factorial series according to (Section 8 of
\cite{We89})
\begin{equation}
s_n \; = \; s \, + \, \omega_n \, 
\sum_{j=0}^{\infty} c_j / (n+\zeta)_j \, ,
\end{equation} 
where $(n+\zeta)_j = \Gamma(n+j+\zeta)/\Gamma(n+\zeta)$ is a Pochhammer
symbol (p.\ 3 of \cite{MaObSo66}). In this way, the sequence
transformation (Eq.\ (8.2-7) of
\cite{We89})
\begin{eqnarray}
\lefteqn{{\cal S}_{k}^{(n)} (\zeta , s_n, \omega_n)}
\nonumber \\
& & \; = \; \frac
{\displaystyle
\sum_{j=0}^{k} \; ( - 1)^{j} \; {{k} \choose {j}} \;
\frac {(\zeta + n +j )_{k-1}} {(\zeta + n + k )_{k-1}} \;
\frac {s_{n+j}} {\omega_{n+j}} }
{\displaystyle
\sum_{j=0}^{k} \; ( - 1)^{j} \; {{k} \choose {j}} \;
\frac {(\zeta + n +j )_{k-1}} {(\zeta + n + k )_{k-1}} \;
\frac {1} {\omega_{n+j}} } 
\label{SidTr}
\end{eqnarray}
can be derived which is formally very similar to Levin's sequence
transformation. The only difference between the transformations ${\cal
L}_{k}^{(n)} (\zeta , s_n, \omega_n)$ and ${\cal S}_{k}^{(n)} (\zeta ,
s_n, \omega_n)$ is that the powers $(\zeta + n +j )^{k-1}$ in
(\ref{LevTr}) are replaced by the Pochhammer symbols $(\zeta + n +
j)_{k-1}$ in (\ref{SidTr}). Again, the shift parameter $\zeta$ has to be
positive in order to admit $n = 0$ in (\ref{SidTr}), and the most
obvious choice is also $\zeta = 1$ which is exclusively used in this
article. Recurrence formulas for the numerator and denominator sums of
${\cal S}_{k}^{(n)} (\zeta, s_n, \omega_n)$ can be found in Section 8.3
of \cite{We89}.

If we use the remainder estimate (\ref{dRemEst}) also in (\ref{SidTr}),
we obtain the following sequence transformation (Eq.\ (8.4-4) of
\cite{We89}):
\begin{equation}
{\delta}_{k}^{(n)} (\zeta, s_n) \; = \; 
{\cal S}_{k}^{(n)} (\zeta, s_n, \Delta s_n) \, .
\label{dSidTr}
\end{equation}
It was shown in several articles that the transformation (\ref{SidTr})
as well as its variant (\ref{dSidTr}) can be very effective
\cite{We89,We90,WeCi90,We94,We96d,JeMoSoWe99,WeSte89,WeCiVi91,We92,%
WeCiVi93,We96a,We96b,We96c,We97,JeBeWeS02000,RoBhaBho96,BhaRoBho97,%
RoBhaBho98,SarSenHalRoy98}, in particular if strongly divergent
alternating series are to be summed.

In the case of the transformations (\ref{dLevTr}) and (\ref{dSidTr}),
the approximation to the limit with the highest transformation order is
given by
\begin{equation}
\{ s_0, s_1, \ldots ,s_{m+1} \} \; \longrightarrow \; 
\Xi_{m}^{(0)} (\zeta, s_0) \, ,
\label{LevTypeApprLim}
\end{equation}
where $\Xi_{k}^{(n)} (\zeta, s_n)$ stands for either $d_{k}^{(n)}
(\zeta, s_n)$ or ${\delta}_{k}^{(n)} (\zeta, s_n)$.

If the input data $s_n$ are the partial sums of a (formal) power series
for some function $f (z)$ according to (\ref{PowSerPS}), $s_n = f_n
(z)$, then the transformations (\ref{dLevTr}) and (\ref{dSidTr}) produce
rational functions $d_{k}^{(n)} \bigr( \zeta, f_n (z) \bigl)$ and
${\delta}_{k}^{(n)} \bigr( \zeta, f_n (z) \bigl)$, whose numerator and
denominator polynomials are of degrees $k+n$ and $k$ in $z$,
respectively (Eqs.\ (4.25) and (4.26) of \cite{WeCiVi93}). Moreover, the
rational approximants $d_{k}^{(n)} \bigr( \zeta, f_n (z) \bigl)$ and
${\delta}_{k}^{(n)} \bigr( \zeta, f_n (z) \bigl)$ satisfy the following
asymptotic error estimates as $z \to 0$ (Eqs.\ (4.28) and (4.29) of
\cite{WeCiVi93}), 
\begin{eqnarray}
f (z) \, - \, d_k^{(n)} \bigl(\zeta, f_n (z) \bigr)
& \; = \; & O (z^{k + n + 2}) \, ,
\label{OdLevPS} \\
f (z) \, - \, {\delta}_k^{(n)} \bigl(\zeta, f_n (z) \bigr)
& \; = \; & O (z^{k + n + 2}) \, ,
\label{OdSidPS}
\end{eqnarray}
which are very similar to the well known accuracy-through-order
relationships satisfied by Pad\'{e} approximants \cite{BaGM96}.

It is a typical feature of all sequence transformations discussed in
this Appendix that they are both {\em homogeneous\/} and {\em
translative}: If the elements of two sequences $\{ s_n
\}_{n=0}^{\infty}$ and $\{ \sigma_n \}_{n=0}^{\infty}$ satisfy
\begin{equation}
\sigma_n \; = \; a s_n \, + \, b \, ,
\end{equation}
where $a$ and $b$ are suitable constants, then
\begin{equation}
{\cal T} (\sigma_n, \sigma_{n+1}, \ldots) \; = \; 
a \, {\cal T} (s_n, s_{n+1}, \ldots) \, + \, b \, .
\label{quasi_lin}
\end{equation}
Sequence transformations ${\cal T}$ satisfying this condition are called 
{\em quasi-linear\/} in the book by Brezinski and Redivo Zaglia
\cite{BreRZa91}. In Section 1.4 of this book, a detailed discussion of 
the properties of quasi-linear sequence transformations as well as
further references can be found.

\section{The Infinite Coupling Limit of the Sextic Anharmonic Oscillator}
\label{App_B}

The detrimental effect of the irregular behavior of the leading elements
of a sequence in convergence acceleration and summation processes is not
restricted to mathematical model problems but occurs also in the
mathematical treatment of scientific problems.

This will be shown by performing extensive summation calculations for
the so-called strong coupling limit $k_3$ of the sextic anharmonic
oscillator. The quartic ($m = 2$), sextic ($m = 3$), and octic ($m = 4$)
anharmonic oscillators are defined by the Hamiltonians
\begin{equation}
\hat{H} (\beta) \; = \; \hat{p}^2 \, + \, \hat{x}^2 \, + \, 
\beta \hat{x}^{2 m} \, , \qquad m = 2, 3, 4 \, ,
\end{equation}
and the strong coupling limit $k_m$ of the ground state energy
eigenvalue $E^{(m} (\beta)$ of this Hamiltonian is defined by
\begin{equation}
k_m \; = \; \lim_{\beta \to \infty} \, 
E^{(m)} (\beta)/\beta^{1/(m+1)} \, .
\label{InfCoupLim}
\end{equation}

Ever since the seminal work of Bender and Wu
\cite{BeWu69,BeWu71,BeWu73}, the divergent weak coupling perturbation
expansion
\begin{equation}
E^{(m)} (\beta) \; = \; 
\sum_{n=0}^{\infty} \, b_{n}^{(m)} \, \beta^n 
\label{wcPT_Em}
\end{equation}
for the ground state energy of an anharmonic oscillator has been
considered to be the model example of a strongly divergent quantum
mechanical perturbation expansion which has to be summed in order to
produce numerically useful results. Accordingly, there is an extensive
literature on the summation of the divergent perturbation
expansions of the anharmonic oscillators (see for example
\cite{WeCi90,We96b,CiWeBraSpi96,We96c,We97,BeWu69,BeWu71,BeWu73,Si70,%
SkaCiZa99a,SkaCiZa99b,ViCi91} and references therein).

In addition to the divergent weak coupling expansion (\ref{wcPT_Em}),
there is also a strong coupling expansion \cite{Si70}
\begin{equation}
E^{(m)} (\beta) \; = \; \beta^{1/(m+1)} \,
\sum_{n=0}^{\infty} \, K_{n}^{(m)} \, \beta^{-2n/(m+1)} \, .
\label{scPT_Em}
\end{equation}
It can be shown that this expansion converges for sufficiently large
values of $\beta$ \cite{Si70,SkaCiZa99a,SkaCiZa99b}. However, the
computation of the perturbative coefficients $K_{n}^{(m)}$ is very
difficult (see for example \cite{We96c} and references therein).

It follows from (\ref{InfCoupLim}), (\ref{wcPT_Em}), and (\ref{scPT_Em})
that the infinite coupling limit $k_m$ corresponds to the leading
coefficient of the strong coupling expansion (\ref{scPT_Em}) according
to
\begin{equation}
k_m \; = \; K_{0}^{(m)} \, .
\end{equation}

The weak coupling perturbation expansion (\ref{wcPT_Em}) cannot be used
in a straightforward way for a calculation of the strong coupling limit
$k_m$. However, this can be accomplished comparatively easily with the
help of the following {\em renormalized\/} weak coupling perturbation
expansion (Eqs.\ (3.30) - (3.31) of \cite{WeCiVi93}):
\begin{equation}
E^{(m)} (\beta) \; = \; (1 - \kappa)^{-1/2} \, 
\sum_{n=0}^{\infty} \, c_{n}^{(m)} \, \kappa^n \, .
\label{rwcEm}
\end{equation}
This expansion is based on a renormalization scheme introduced by
Vinette and \v{C}\'{\i}\v{z}ek \cite{ViCi91}. In this approach, the
original coupling constant $\beta \in [0, \infty)$ is transformed into a
renormalized and explicitly $m$-dependent coupling constant $\kappa \in
[0, 1)$ according to (Eq.\ (3.19) of \cite{WeCiVi93})
\begin{equation}
\beta \; = \; \frac {1} { B_m } \,
\frac {\kappa} {(1 - \kappa)^{(m + 1)/2}} \, ,
\end{equation}
where (Eq.\ (3.17) of \cite{WeCiVi93})
\begin{equation}
B_m \; = \; \frac {m \, (2 m - 1)!!} {2^{m - 1}} \, .
\end{equation}
For the sextic ($m = 3$) case, these expressions correspond to $B_3 =
45/4$ and $\beta = 4 \kappa/[45(1-\kappa)^2]$.

Thus, the infinite coupling limits $k_3$ of the sextic anharmonic
oscillator can be expressed by the renormalized weak coupling expansion
(\ref{rwcEm}) according to (Eqs.\ (3.43) and (3.44) of \cite{WeCiVi93})
\begin{equation}
k_3 \; = \; [45/4]^{1/4} \, \sum_{n=0}^{\infty} \, c_{n}^{(3)} \, .
\label{rPT_k3}
\end{equation}
The summation of either this or the perturbation series (\ref{rwcEm})
with $m = 3$, from which (\ref{rPT_k3}) was derived, is a formidable
computational problem. This follows at once from the large-$n$
asymptotics of the renormalized perturbative coefficients for the sextic
anharmonic oscillator (Eq.\ (3.34) of \cite{WeCiVi93}):
\begin{eqnarray}
\lefteqn{c_{n}^{(3)} \, \sim \, (-1)^{n+1} \,
\frac{(128)^{1/2}}{\pi^{2}}} \nonumber \\
& & \quad \times \, \Gamma(2 n + 1/2) \,
\bigl( 64 / [45 \pi^{2}] \bigr)^{n} \, , \qquad n \to \infty \, .
\label{cAs3}
\end{eqnarray}
It should be noted that the summation of the perturbation series
(\ref{rPT_k3}) for $k_3$ is much more demanding than the summation of
the divergent asymptotic expansions for special functions since their
coefficients $c_n$ grow essentially like $n!$
\cite{We89,We90,WeCi90,We96d}.  Although Pad\'e approximants -- or
equivalently Wynn's epsilon algorithm -- are in principle capable of
summing alternating divergent power series whose coefficients $c_n$ grow
essentially like $(2 n)!$ in magnitude, the convergence of Pad\'e
approximants is too slow to be practically useful. Moreover, it was
shown in \cite{WeCiVi93} that the Levin transformation (\ref{dLevTr})
produces in the case of the perturbation expansions for the anharmonic
oscillators sequences of approximants which initially seem to converge
but which ultimately diverge. In contrast, the Levin-type transformation
(\ref{dSidTr}) produces comparatively good results.

Thus, in analogy to \cite{WeCiVi93} we sum the perturbation series
(\ref{rPT_k3}) with the help of the Levin-type transformation
${\delta}_{k}^{(n)} (\zeta, s_n)$ defined in (\ref{dSidTr}). It was
shown in several articles that the sequence transformation
(\ref{dSidTr}) as well as the transformation (\ref{SidTr}), from which
it was derived, is apparently very effective, in particular if strongly
divergent alternating series are to be summed
\cite{We89,We90,WeCi90,We94,We96d,JeMoSoWe99,WeSte89,WeCiVi91,We92,%
WeCiVi93,We96a,We96b,We96c,We97,JeBeWeS02000,RoBhaBho96,BhaRoBho97,%
RoBhaBho98,SarSenHalRoy98}.

In our summation calculations for $k_3$ we use the renormalized
coefficients $c_{\nu}^{(3)}$ with $0 \le \nu \le 300$. The coefficients
were calculated using the exact rational arithmetics of Maple by solving
a system of nonlinear difference equations as described in the Appendix
of \cite{WeCiVi93}. Unfortunately, Eq.\ (A22) in the Appendix of
\cite{WeCiVi93}, which specifies the system of nonlinear equations,
contains a typographical error. Correct is
\begin{eqnarray}
4 j \, G_j^{(n)} \; & = \; &
2 (j + 1) (2 j + 1) \, G_{j+1}^{(n)} \, + \,
\frac {1} {B_m} G_{j - m}^{(n-1)} \nonumber \\
& & - \,  G_{j - 1}^{(n-1)} \, - \,
2 \, \sum_{k=1}^{n-1} \, G_1^{(k)} \, G_j^{(n-k)} \, .
\end{eqnarray}

The topic of this article is the study of the impact of irregular input
data on the performance of sequence transformations. Accordingly, we
have to investigate whether the renormalized coefficients $c_{n}^{(3)}$
behave irregularly for small indices $n$. For that purpose, we list in
Table \ref{Tab_b_1} selected renormalized coefficients $c_{n}^{(3)}$ as
well as the corresponding ratios
\begin{equation}
{\cal C}_{n}^{(3)} \; = \; 
\frac{(-1)^{n+1} \, \pi^{2} \, c_{n}^{(3)}}
{\sqrt{128} \, \Gamma(2 n + 1/2)} \,
\left(\frac{45 \pi^{2}}{64}\right)^n \, ,
\label{Asy_n_C3}
\end{equation}
that are obtained by dividing the renormalized coefficients
$c_{n}^{(3)}$ by the leading order of their large-$n$ asymptotics
according to (\ref{cAs3}).

The last column in Table \ref{Tab_b_1} shows quite clearly that the
renormalized coefficients $c_{n}^{(3)}$ deviate for small values of $n$
considerably from their larger-order behavior. Firstly, the coefficients
$c_{0}^{(3)}$ and $c_{1}^{(3)}$ apparently possess ``wrong''
sign. Secondly, the coefficients $c_{n}^{(3)}$ initially decrease in
magnitude, and only for $n \ge 4$ they grow as they should according to
(\ref{cAs3}). Nevertheless, it is only a relatively mild irregularity,
which affects only a few of the available coefficients $c_{\nu}^{(3)}$
with $0 \le \nu \le 300$.

The impact of the irregular coefficients $c_{n}^{(3)}$ with small
indices $n$ can be checked by computing for $l = 0, 1, 2, \ldots$ and
for $n \le 299 - l$ the approximants
\begin{equation}
k_3^{(n, l)} \; = \; \delta_{n}^{(0)} \bigl(1, s_{0}^{(l)} \bigr) \, ,
\label{k3_ln}
\end{equation} 
to the infinite coupling limit $k_3$ of the sextic anharmonic
oscillator. Here, $\delta_{n}^{(0)}$ is the Levin-type transformation
(\ref{dSidTr}), and
\begin{equation}
s_{n}^{(l)} \; = \; 
[45/4]^{1/4} \, \sum_{\nu=0}^{n+l} \, c_{\nu}^{(3)}
\end{equation} 
is a partial sum of the perturbation series (\ref{rPT_k3}) which skips
the first $l$ terms in the transformation process.

In Table \ref{Tab_b_2}, the approximants $k_3^{(l, n)}$ with the three
highest possible values of $n \le 299 - l$ are listed for $l \le 12$. 

If we compare the results with the extremely accurate result of Vinette
and {\v{C}\'{\i}\v{z}ek} (Eq.\ (69) of \cite{ViCi91})
\begin{equation}
 k_3 \; = \; 1.144~802~453~797~052~763~765~457~534~149~549 \, ,
\end{equation}
which was obtained nonperturbatively, we see that we gain 5 decimal
digits by skipping the first 7 terms of the perturbation series
(\ref{rPT_k3}) for $k_3$. For $l \ge 8$, the accuracy of the summation
results deteriorates again.

This is a very remarkable gain of accuracy, if we take into account that
the summation of the perturbation series (\ref{rPT_k3}) is a formidable
task and that the leading terms $c_{n}^{(3)}$ display an only relatively
mild irregularity, as shown in Table \ref{Tab_b_1}. Nevertheless, the
results show quite clearly that the transformation order is not the only
criterion which affects the performance of a sequence
transformation. The results in Table \ref{Tab_b_2} show that may be more
effective to use smaller sets of more regular input data.

\section{Recurrence Formulas for the Gaussian Hypergeometric Function} 
\label{App_C}

Many three-term recurrence formulas satisfied by the Gaussian
hypergeometric function ${}_2 F_1 (a, b; c; z)$ are known. For example,
on pp.\ 557 - 558 of \cite{AbSte72} or on pp.\ 46 - 47 of
\cite{MaObSo66}, the following formulas can be found:
\begin{eqnarray}
\lefteqn{(c-a) \, {}_2 F_1 (a-1, b; c; z)} \nonumber \\
& + & \, [2a-c-(a-b)z] \, {}_2 F_1 (a, b; c; z) \nonumber \\ 
& & + \; a(z-1) \, {}_2 F_1 (a+1, b; c; z) \; = \; 0 \, ,
\\
\lefteqn{(c-b) \, {}_2 F_1 (a, b-1; c; z)} \nonumber \\
& + & [2b-c-(b-a)z] \, {}_2 F_1 (a, b; c; z) \nonumber \\ 
& & + \; b(z-1) \, {}_2 F_1 (a, b+1; c; z) \; = \; 0 \, ,
\\
\lefteqn{c(1-c)(1-z) \, {}_2 F_1 (a, b; c-1; z)} \nonumber \\
& + & \, c[c-1-(2c-a-b-1)z] \, {}_2 F_1 (a, b; c; z) \nonumber \\
& & + \; (c-a)(c-b)z \, {}_2 F_1 (a, b; c+1; z) \; = \; 0 \, ,
\label{Rec_2F1_3}
\\
\lefteqn{(b-a) \, {}_2 F_1 (a, b; c; z)
\, + \, a \, {}_2 F_1 (a+1, b; c; z)} \nonumber \\
& & - \; b \, {}_2 F_1 (a, b+1; c; z) \; = \; 0 \, ,
\\
\lefteqn{(b-c) \, {}_2 F_1 (a, b-1; c; z)} \nonumber \\
& + & (c-a-b) \, {}_2 F_1 (a, b; c; z) \nonumber \\
& & + \; a(1-z) \, {}_2 F_1 (a+1, b; c; z) \; = \; 0 \, ,
\\
\lefteqn{c[a-(c-b)z] \, {}_2 F_1 (a, b; c; z)} \nonumber \\
& - & \, ac(1-z) \, {}_2 F_1 (a+1, b; c; z) \nonumber \\
& & + \; (c-a)(c-b)z \, {}_2 F_1 (a, b; c+1; z) \; = \; 0 \, ,
\label{Rec_2F1_6}
\\
\lefteqn{(1-c) \, {}_2 F_1 (a, b; c-1; z)} \nonumber \\
& + & \, (c-a-1) \, {}_2 F_1 (a, b; c; z) \nonumber \\
& & + \; a \, {}_2 F_1 (a+1, b; c; z) \; = \; 0 \, ,
\label{Rec_2F1_7}
\\
\lefteqn{(a-c) \, {}_2 F_1 (a-1, b; c; z)} \nonumber \\
& + & (c-a-b) \, {}_2 F_1 (a, b; c; z) \nonumber \\
& & + \; b(1-z) \, {}_2 F_1 (a, b+1; c; z) \; = \; 0 \, ,
\\
\lefteqn{(a-c) \, {}_2 F_1 (a-1, b; c; z)} \nonumber \\
& + & (c-b) \, {}_2 F_1 (a, b-1; c; z) \nonumber \\
& & + \; (b-a)(1-z) \, {}_2 F_1 (a, b; c; z) \; = \; 0 \, ,
\\
\lefteqn{(-c) \, {}_2 F_1 (a-1, b; c; z)} \nonumber \\
& + & \, c(1-z) \, {}_2 F_1 (a, b; c; z) \nonumber \\
& & + \; (c-b)z \, {}_2 F_1 (a, b; c+1; z) \; = \; 0 \, ,
\label{Rec_2F1_10}
\\
\lefteqn{(c-a) \, {}_2 F_1 (a-1, b; c; z)} \nonumber \\
& - & \, (c-1)(1-z) \, {}_2 F_1 (a, b; c-1; z) \nonumber \\
& & + \; [a-1-(c-b-1)z] \, {}_2 F_1 (a, b; c; z) \; = \; 0 \, ,
\label{Rec_2F1_11}
\\
\lefteqn{c[b-(c-a)z] \, {}_2 F_1 (a, b; c; z)} \nonumber \\
& - & \, bc(1-z) \, {}_2 F_1 (a, b+1; c; z) \nonumber \\
& & + \; (c-a)(c-b)z \, {}_2 F_1 (a, b; c+1; z) \; = \; 0 \, ,
\label{Rec_2F1_12}
\\
\lefteqn{(1-c) \, {}_2 F_1 (a, b; c-1; z)} \nonumber \\
& + & \, (c-b-1) \, {}_2 F_1 (a, b; c; z) \nonumber \\
& & + \; b \, {}_2 F_1 (a, b+1; c; z) \; = \; 0 \, ,
\\
\lefteqn{(-c) \, {}_2 F_1 (a, b-1; c; z)} \nonumber \\
& + & \, c(1-z) \, {}_2 F_1 (a, b; c; z) \nonumber \\
& & + \; (c-a)z \, {}_2 F_1 (a, b; c+1; z) \; = \; 0 \, ,
\label{Rec_2F1_14}
\\
\lefteqn{(c-b) \, {}_2 F_1 (a, b-1; c; z)} \nonumber \\
& - & \, (c-1)(1-z) \, {}_2 F_1 (a, b; c-1; z) \nonumber \\
& & + \; [b-1-(c-a-1)z] \, {}_2 F_1 (a, b; c; z) \; = \; 0 \, .
\label{Rec_2F1_15}
\end{eqnarray}
It is a typical feature of these recurrence formulas that there is a
hypergeometric function ${}_2 F_1 (a, b; c; z)$ plus two other ${}_2
F_1$'s which differ with respect to only {\em one} of the three
parameters by $\pm 1$. However, recurrence formulas, which contain a
hypergeometric function ${}_2 F_1 (a, b; c; z)$ plus two other ${}_2
F_1$'s which differ with respect to {\em two} or even {\em three}
parameters by $\pm 1$, can be constructed comparatively easily. For that
purpose, we combine those of the recurrence formulas given above, whose
the third parameter $c$ assumes at least two different values, with the
linear transformation formulas (see for example p.\ 559 of
\cite{AbSte72} or p.\ 47 of \cite{MaObSo66})
\begin{eqnarray}
{}_2 F_1 (a, b; c; z) & = &
(1-z)^{c-a-b} \, {}_2 F_1 (c-a, c-b; c; z)
\label{LTr_0}
\\
& = &
(1-z)^{-a} \, {}_2 F_1 \bigl( a, c-b; c; z/(z-1) \bigr)
\label{LTr_1}
\\
& = &
(1-z)^{-b} \, {}_2 F_1 \bigl( c-a, b; c; z/(z-1) \bigr) \, .
\label{LTr_2}
\end{eqnarray}

For the derivation of new recurrence formulas, we replace in
(\ref{Rec_2F1_3}) $a$ by $c-a$ and $b$ by $c-b$. This yields:
\begin{eqnarray}
\lefteqn{c(1-c)(1-z) \, {}_2 F_1 (c-a, c-b; c-1; z)} \nonumber \\
& + & \, c[c-1-(a+b-1) z] \, {}_2 F_1 (c-a, c-b; c; z) \nonumber \\
& & + \; 
a b z \, {}_2 F_1 (c-a, c-b; c+1; z) \; = \; 0 \, .
\end{eqnarray}
If we now combine this relationship with the linear transformation
(\ref{LTr_0}), we obtain the following recurrence formula, where all
three parameters of the hypergeometric functions change simultaneously: 
\begin{eqnarray}
\lefteqn{c(1-c) \, {}_2 F_1 (a-1, b-1; c-1; z)} \nonumber \\
& + & \, c[c-1-(a+b-1) z] \, {}_2 F_1 (a, b; c; z) \nonumber \\
& & + \; 
a b z (1-z) \, {}_2 F_1 (a+1, b+1; c+1; z) \; = \; 0 \, .
\end{eqnarray}
If we now proceed in (\ref{Rec_2F1_6}), (\ref{Rec_2F1_7}), and in
(\ref{Rec_2F1_10}) - (\ref{Rec_2F1_15}) in exactly the same way, we
obtain the following recurrence formulas:
\begin{eqnarray}
\lefteqn{c(c-a-bz) \, {}_2 F_1 (a, b; c; z)} \nonumber \\
& - & \, c(c-a) \, {}_2 F_1 (a-1, b; c; z) \nonumber \\
& & + \; abz(1-z) \, {}_2 F_1 (a+1, b+1; c+1; z) \; = \; 0 \, ,
\label{Rec_abc2}
\\
\lefteqn{(1-c) \, {}_2 F_1 (a-1, b-1; c-1; z)} \nonumber \\
& + & \, (a-1)(1-z) \, {}_2 F_1 (a, b; c; z) \nonumber \\
& & + \; (c-a) \, {}_2 F_1 (a-1, b; c; z) \; = \; 0 \, ,
\label{Rec_abc1}
\\
\lefteqn{c \, {}_2 F_1 (a+1, b; c; z)} \nonumber \\
& - & \, c \, {}_2 F_1 (a, b; c; z) \nonumber \\
& & - \; bz \, {}_2 F_1 (a+1, b+1; c+1; z) \; = \; 0 \, ,
\label{Rec_abc3}
\\
\lefteqn{a(1-z) \, {}_2 F_1 (a+1, b; c; z)} \nonumber \\
& + & \, (1-c) \, {}_2 F_1 (a-1, b-1; c-1; z) \nonumber \\
& & + \; [c-a-1-(b-1)z] \, {}_2 F_1 (a, b; c; z) \; = \; 0 \, ,
\label{Rec_abc4}
\\
\lefteqn{c(c-b-az) \, {}_2 F_1 (a, b; c; z)} \nonumber \\
& - & \, c(c-b) \, {}_2 F_1 (a, b-1; c; z) \nonumber \\
& & + \; a b z (1-z) \, {}_2 F_1 (a+1, b+1; c+1; z) \; = \; 0 \, ,
\label{Rec_abc5}
\\
\lefteqn{(1-c) \, {}_2 F_1 (a-1, b-1; c-1; z)} \nonumber \\
& + & \, (b-1)(1-z) \, {}_2 F_1 (a, b; c; z) \nonumber \\
& & + \; (c-b) \, {}_2 F_1 (a, b-1; c; z) \; = \; 0 \, ,
\label{Rec_abc6}
\\
\lefteqn{c \, {}_2 F_1 (a, b+1; c; z)} \nonumber \\
& - & \, c \, {}_2 F_1 (a, b; c; z) \nonumber \\
& & - \; az \, {}_2 F_1 (a+1, b+1; c+1; z) \; = \; 0 \, ,
\label{Rec_abc7}
\\
\lefteqn{b(1-z) \, {}_2 F_1 (a, b+1; c; z)} \nonumber \\
& + & \, (1-c) \, {}_2 F_1 (a-1, b-1; c-1; z) \nonumber \\
& & + \; [c-b-1-(a-1)z] \, {}_2 F_1 (a, b; c; z) \; = \; 0 \, .
\label{Rec_abc8}
\end{eqnarray}
Not all of these recurrence formulas are independent. For example,
(\ref{Rec_abc2}) and (\ref{Rec_abc5}) can be transformed into each other
by interchanging $a$ and $b$. This is also true for the pairs
(\ref{Rec_abc1}) and (\ref{Rec_abc6}), (\ref{Rec_abc3}) and
(\ref{Rec_abc7}), and (\ref{Rec_abc4}) and (\ref{Rec_abc8}), which can
be transformed into each other by interchanging $a$ and $b$.

For the derivation of recurrence formulas which differ with respect to
two parameters by $\pm 1$, we replace in (\ref{Rec_2F1_3}) $a$ by $c-a$
and $z$ by $z/(z-1)$. This yields:
\begin{eqnarray}
\lefteqn{\frac{c(1-c)}{1-z} \, 
{}_2 F_1 \left(c-a, b; c-1; \frac{z}{z-1}\right)} \nonumber \\
& + & \, \frac{c[c-1+(a-b)z]}{z-1} \, 
{}_2 F_1 \left(c-a, b; c; \frac{z}{z-1}\right) \nonumber \\
& & + \; \frac{a(c-b)z}{z-1} \, 
{}_2 F_1 \left(c-a, b; c+1; \frac{z}{z-1}\right) \; = \; 0 \, .
\end{eqnarray}
If we now combine this relationship with the linear transformation
(\ref{LTr_2}), we obtain the following recurrence formula, where the
first and the third parameter of the hypergeometric series change
simultaneously:
\begin{eqnarray}
\lefteqn{c(1-c) \, {}_2 F_1 (a-1, b; c-1; z)} \nonumber \\
& - & \, c[c-1+(a-b)z] \, {}_2 F_1 (a, b; c; z) \nonumber \\
& & + \; 
a(c-b)z \, {}_2 F_1 (a+1, b; c+1; z) \; = \; 0 \, .
\label{Rec_ac}
\end{eqnarray}
If we now proceed in (\ref{Rec_2F1_6}), (\ref{Rec_2F1_7}), and in
(\ref{Rec_2F1_10}) - (\ref{Rec_2F1_15}) in exactly the same way, we
obtain the following recurrence formulas:
\begin{eqnarray}
\lefteqn{c[c-a+(a-b)z] \, {}_2 F_1 (a, b; c; z)} \nonumber \\
& - & \, c(c-a) \, {}_2 F_1 (a-1, b; c; z) \nonumber \\
& & - \; a(c-b)z \, {}_2 F_1 (a+1, b; c+1; z) \; = \; 0 \, ,
\label{Rec_ac2}
\\
\lefteqn{(1-c) \, {}_2 F_1 (a-1, b; c-1; z)} \nonumber \\
& + & \, (a-1) \, {}_2 F_1 (a, b; c; z) \nonumber \\
& & + \; (c-a) \, {}_2 F_1 (a-1, b; c; z) \; = \; 0 \, ,
\label{Rec_ac1}
\\
\lefteqn{c(1-z) \, {}_2 F_1 (a+1, b; c; z)} \nonumber \\
& - & \, c \, {}_2 F_1 (a, b; c; z) \nonumber \\
& & + \; (c-b)z \, {}_2 F_1 (a+1, b; c+1; z) \; = \; 0 \, ,
\label{Rec_ac3}
\\
\lefteqn{a(1-z) \, {}_2 F_1 (a+1, b; c; z)} \nonumber \\
& + & \, (1-c) \, {}_2 F_1 (a-1, b; c-1; z) \nonumber \\
& & + \; [c-a-1+(a-b)z] \, {}_2 F_1 (a, b; c; z) \; = \; 0 \, ,
\label{Rec_ac4}
\\
\lefteqn{c[b+(a-b)z] \, {}_2 F_1 (a, b; c; z)} \nonumber \\
& - & \, bc(1-z) \, {}_2 F_1 (a, b+1; c; z) \nonumber \\
& & - \; a(c-b)z \, {}_2 F_1 (a+1, b; c+1; z) \; = \; 0 \, ,
\label{Rec_ac5}
\\
\lefteqn{(1-c) \, {}_2 F_1 (a-1, b; c-1; z)} \nonumber \\
& + & \, (c-b-1) \, {}_2 F_1 (a, b; c; z) \nonumber \\
& & + \; b(1-z) \, {}_2 F_1 (a, b+1; c; z) \; = \; 0 \, ,
\label{Rec_ac6}
\\
\lefteqn{c \, {}_2 F_1 (a, b-1; c; z)} \nonumber \\
& - & \, c \, {}_2 F_1 (a, b; c; z) \nonumber \\
& & + \; az \, {}_2 F_1 (a+1, b; c+1; z) \; = \; 0 \, ,
\label{Rec_ac7}
\\
\lefteqn{(c-b) \, {}_2 F_1 (a, b-1; c; z)} \nonumber \\
& + & \, (1-c) \, {}_2 F_1 (a-1, b; c-1; z) \nonumber \\
& & + \; [b-1+(a-b)z] \, {}_2 F_1 (a, b; c; z) \; = \; 0 \, .
\label{Rec_ac8}
\end{eqnarray}

Next, we replace in (\ref{Rec_2F1_3}) $b$ by $c-b$
and $z$ by $z/(z-1)$. This yields:
\begin{eqnarray}
\lefteqn{\frac{c(1-c)}{1-z} \, 
{}_2 F_1 \left(a, c-b; c-1; \frac{z}{z-1}\right)} \nonumber \\
& + & \, \frac{c[c-1-(a-b)z]}{z-1} \, 
{}_2 F_1 \left(a, c-b; c; \frac{z}{z-1}\right) \nonumber \\
& & + \; \frac{(c-a)bz}{z-1} \, 
{}_2 F_1 \left(a, c-b; c+1; \frac{z}{z-1}\right) \; = \; 0 \, .
\end{eqnarray}
If we now combine this relationship with the linear transformation
(\ref{LTr_1}), we obtain the following recurrence formula, where the
second and the third parameter of the hypergeometric series change
simultaneously:
\begin{eqnarray}
\lefteqn{c(1-c) \, {}_2 F_1 (a, b-1; c-1; z)} \nonumber \\
& + & \, c[c-1-(a-b)z] \, {}_2 F_1 (a, b; c; z) \nonumber \\
& & - \; 
(c-a)bz \, {}_2 F_1 (a, b+1; c+1; z) \; = \; 0 \, .
\label{Rec_bc}
\end{eqnarray}
If we now proceed in (\ref{Rec_2F1_6}), (\ref{Rec_2F1_7}), and in
(\ref{Rec_2F1_10}) - (\ref{Rec_2F1_15}) in exactly the same way, we
obtain the following recurrence formulas:
\begin{eqnarray}
\lefteqn{c[(a-b)z-a] \, {}_2 F_1 (a, b; c; z)} \nonumber \\
& + & \, ac(1-z) \, {}_2 F_1 (a+1, b; c; z) \nonumber \\
& & + \; (c-a)bz \, {}_2 F_1 (a, b+1; c+1; z) \; = \; 0 \, ,
\label{Rec_bc2}
\\
\lefteqn{(1-c) \, {}_2 F_1 (a, b-1; c-1; z)} \nonumber \\
& + & \, (c-a-1) \, {}_2 F_1 (a, b; c; z) \nonumber \\
& & + \; a(1-z) \, {}_2 F_1 (a+1, b; c; z) \; = \; 0 \, ,
\label{Rec_bc1}
\\
\lefteqn{c \, {}_2 F_1 (a-1, b; c; z)} \nonumber \\
& - & \, c \, {}_2 F_1 (a, b; c; z) \nonumber \\
& & + \; bz \, {}_2 F_1 (a, b+1; c+1; z) \; = \; 0 \, ,
\label{Rec_bc3}
\\
\lefteqn{(c-a) \, {}_2 F_1 (a-1, b; c; z)} \nonumber \\
& + & \, (1-c) \, {}_2 F_1 (a, b-1; c-1; z) \nonumber \\
& & + \; [a-1-(a-b)z] \, {}_2 F_1 (a, b; c; z) \; = \; 0 \, ,
\label{Rec_bc4}
\\
\lefteqn{c[c-b-(a-b)z] \, {}_2 F_1 (a, b; c; z)} \nonumber \\
& - & \, c(c-b) \, {}_2 F_1 (a, b-1; c; z) \nonumber \\
& & - \; (c-a)bz \, {}_2 F_1 (a, b+1; c+1; z) \; = \; 0 \, ,
\label{Rec_bc5}
\\
\lefteqn{(1-c) \, {}_2 F_1 (a, b-1; c-1; z)} \nonumber \\
& + & \, (b-1) \, {}_2 F_1 (a, b; c; z) \nonumber \\
& & + \; (c-b) \, {}_2 F_1 (a, b-1; c; z) \; = \; 0 \, ,
\label{Rec_bc6}
\\
\lefteqn{c(1-z) \, {}_2 F_1 (a, b+1; c; z)} \nonumber \\
& - & \, c \, {}_2 F_1 (a, b; c; z) \nonumber \\
& & + \; (c-a)z \, {}_2 F_1 (a, b+1; c+1; z) \; = \; 0 \, ,
\label{Rec_bc7}
\\
\lefteqn{b(1-z) \, {}_2 F_1 (a, b+1; c; z)} \nonumber \\
& + & \, (1-c) \, {}_2 F_1 (a, b-1; c-1; z) \nonumber \\
& & + \; [c-b-1-(a-b)z] \, {}_2 F_1 (a, b; c; z) \; = \; 0 \, .
\label{Rec_bc8}
\end{eqnarray}
The two groups of recursions (\ref{Rec_ac}) - (\ref{Rec_ac8}) and
(\ref{Rec_bc}) - (\ref{Rec_bc8}), respectively, are not
independent. They can be transformed into each other by interchanging
$a$ and $b$.

\newpage

\onecolumn
\centerline{\bf \large TABLES}

\widetext

\begin{table}
\caption{Summation of the divergent hypergeometric series
$z {}_2 F_1 (1, 1; 2; -z) = \ln (1+z)$ for $z = 7/2$ with the help of
Wynn's epsilon algorithm, Brezinski's theta algorithm, and the
Levin-type transformations $d_{k}^{(n)} (\zeta, s_n)$ and
$\delta_{k}^{(n)} (\zeta, s_n)$.}
\label{Tab_3_1}
\begin{tabular}{lrrrrr}%
$n$%
& \multicolumn{1}{c}{$s_{n} (z)$}%
& \multicolumn{1}{c}{$\epsilon_{2 \Ent {n/2}}^{(n - 2 \Ent {n/2})}$}%
& \multicolumn{1}{c}{$\theta_{2 \Ent {n/3}}^{(n - 3 \Ent {n/3})}$}%
& \multicolumn{1}{c}{$d_{n}^{(0)} \bigl(1, s_0 (z) \bigr)$}%
& \multicolumn{1}{c}{${\delta}_{n}^{(0)} \bigl(1, s_0 (z) \bigr)$}%
\rule[-6pt]{0pt}{18pt} \\
& \multicolumn{1}{c}{Eq.\ (\protect\ref{ParSumT3_1})}%
& \multicolumn{1}{c}{Eq.\ (\protect\ref{eps_al})}%
& \multicolumn{1}{c}{Eq.\ (\protect\ref{theta_al})}%
& \multicolumn{1}{c}{Eq.\ (\protect\ref{dLevTr})}%
& \multicolumn{1}{c}{Eq.\ (\protect\ref{dSidTr})}%
\rule[-6pt]{0pt}{12pt} 
\\
\hline%
$0 $ & $ 0.350~000 \times 10^{1}$ & $ 3.500~000~000~000$ & $
3.500~000~000~000$ & $3.500~000~000~000$ & $3.500~000~000~000$%
\rule[-1pt]{0pt}{12pt} \\
$1 $ & $-0.262~500 \times 10^{1}$ & $-2.625~000~000~000$ & 
$-2.625~000~000~000$ & $1.662~500~000~000$ & $1.662~500~000~000$ \\
$2 $ & $ 0.116~667 \times 10^{2}$ & $ 1.662~500~000~000$ & 
$11.666~666~666~667$ & $1.471~337~579~618$ & $1.471~337~579~618$ \\
$3 $ & $-0.258~490 \times 10^{2}$ & $ 1.317~528~735~632$ & 
$ 1.561~447~811~448$ & $1.507~573~834~307$ & $1.502~377~638~599$ \\
$4 $ & $ 0.791~948 \times 10^{2}$ & $ 1.521~596~244~131$ & 
$ 1.447~356~630~824$ & $1.504~123~629~505$ & $1.504~105~974~245$ \\
$5 $ & $-0.227~183 \times 10^{3}$ & $ 1.488~926~130~389$ & 
$ 1.579~131~944~444$ & $1.504~012~642~929$ & $1.504~105~974~245$ \\
$6 $ & $ 0.691~950 \times 10^{3}$ & $ 1.506~184~895~833$ & 
$ 1.505~133~549~503$ & $1.504~083~039~403$ & $1.504~083~649~440$ \\
$7 $ & $-0.212~289 \times 10^{4}$ & $ 1.502~565~692~760$ & 
$ 1.503~300~027~415$ & $1.504~078~127~161$ & $1.504~078~473~952$ \\
$8 $ & $ 0.663~440 \times 10^{4}$ & $ 1.504~338~769~814$ & 
$ 1.504~815~718~579$ & $1.504~077~213~984$ & $1.504~077~564~801$ \\
$9 $ & $-0.209~511 \times 10^{5}$ & $ 1.503~911~593~667$ & 
$ 1.504~085~412~192$ & $1.504~077~398~663$ & $1.504~077~421~623$ \\
$10$ & $ 0.668~209 \times 10^{5}$ & $ 1.504~110~270~747$ & 
$ 1.504~071~710~932$ & $1.504~077~400~367$ & $1.504~077~400~332$ \\
$11$ & $-0.214~781 \times 10^{6}$ & $ 1.504~058~264~517$ & 
$ 1.504~082~144~389$ & $1.504~077~396~419$ & $1.504~077~397~274$ \\
$12$ & $ 0.695~009 \times 10^{6}$ & $ 1.504~081~563~208$ & 
$ 1.504~077~411~492$ & $1.504~077~396~733$ & $1.504~077~396~845$ \\
$13$ & $-0.226~181 \times 10^{7}$ & $ 1.504~075~121~301$ & 
$ 1.504~077~381~028$ & $1.504~077~396~787$ & $1.504~077~396~786$ \\
$14$ & $ 0.739~713 \times 10^{7}$ & $ 1.504~077~927~290$ & 
$ 1.504~077~411~335$ & $1.504~077~396~776$ & $1.504~077~396~778$ \\
$15$ & $-0.242~963 \times 10^{8}$ & $ 1.504~077~120~833$ & 
$ 1.504~077~396~437$ & $1.504~077~396~776$ & $1.504~077~396~776$ \\
\hline
Exact & & $1.504~077~396~776$ & $1.504~077~396~776$ & 
$1.504~077~396~776$ & $1.504~077~396~776$%
\rule[-1pt]{0pt}{12pt} \\
\end{tabular}
\end{table}

\begin{table}
\caption{Summation of the hypergeometric series 
${}_2 F_1 (2/3, 4/3; 1/3; z) = (1+z) (1-z)^{-5/3}$ for $z = (1+{\rm
i}\sqrt{3})/2$ with the help of Wynn's epsilon algorithm and Levin's
transformation $d_{k}^{(n)} (\zeta, s_n)$.}
\label{Tab_3_2}
\begin{tabular}{lrrr}%
$n$ & \multicolumn{1}{c}{$s_{n} (z)$}
& \multicolumn{1}{c}{$\epsilon_{2 \Ent {n/2}}^{(n - 2 \Ent {n/2})}$}%
& \multicolumn{1}{c}{$d_{n}^{(0)} \bigl(1, s_0 (z) \bigr)$}%
\rule[-6pt]{0pt}{18pt} \\
 & \multicolumn{1}{c}{Eq.\ (\protect\ref{ParSumT3_2})} 
& \multicolumn{1}{c}{Eq.\ (\protect\ref{eps_al})}
& \multicolumn{1}{c}{Eq.\ (\protect\ref{dLevTr})} \rule[-6pt]{0pt}{12pt} \\
\hline%
$0 $ & \multicolumn{1}{l}{$\phantom{-}$ $\!1.00~000$} &
\multicolumn{1}{l}{$\phantom{-}$ $\!1.000~000~000$} &
\multicolumn{1}{l}{$\phantom{-}$ $\!1.000~000~000~000$%
\rule[-1pt]{0pt}{12pt}} \\ 
$1 $ & $ 2.33~333 \, + \, {\rm i} \, 2.30~940$ & 
$ 2.333~333~334 \, + \, {\rm i} \, 2.309~401~077$ & 
$-0.531~737~773~153 \, + \, {\rm i} \, 1.384~198~772~334$ \\      
$2 $ & $ 0.38~889 \, + \, {\rm i} \, 5.67~728$ & 
$-0.531~737~773 \, + \, {\rm i} \, 1.384~198~772$ & 
$-1.131~921~176~770 \, + \, {\rm i} \, 1.411~935~517~843$ \\      
$3 $ & $-4.54~938 \, + \, {\rm i} \, 5.67~728$ & 
$-0.953~868~768 \, + \, {\rm i} \, 1.632~539~413$ & 
$-1.120~276~444~364 \, + \, {\rm i} \, 1.322~600~394~819$ \\      
$4 $ & $-7.49~177 \, + \, {\rm i} \, 0.58~091$ & 
$-1.136~352~936 \, + \, {\rm i} \, 1.326~189~971$ & 
$-1.112~622~137~231 \, + \, {\rm i} \, 1.326~484~602~288$ \\      
$5 $ & $-4.11~180 \, - \, {\rm i} \, 5.27~337$ & 
$-1.115~979~051 \, + \, {\rm i} \, 1.318~232~359$ & 
$-1.113~364~776~960 \, + \, {\rm i} \, 1.326~919~577~184$ \\      
$6 $ & $ 3.46~967 \, - \, {\rm i} \, 5.27~337$ & 
$-1.111~974~302 \, + \, {\rm i} \, 1.326~730~914$ & 
$-1.113~348~742~255 \, + \, {\rm i} \, 1.326~817~605~791$ \\      
$7 $ & $ 7.64~993 \, + \, {\rm i} \, 1.96~705$ & 
$-1.113~209~082 \, + \, {\rm i} \, 1.327~290~958$ & 
$-1.113~338~972~201 \, + \, {\rm i} \, 1.326~827~801~234$ \\      
$8 $ & $ 3.09~756 \, + \, {\rm i} \, 9.85~197$ & 
$-1.113~429~777 \, + \, {\rm i} \, 1.326~838~273$ & 
$-1.113~340~931~714 \, + \, {\rm i} \, 1.326~828~116~649$ \\      
$9 $ & $-6.72~205 \, + \, {\rm i} \, 9.85~197$ & 
$-1.113~348~679 \, + \, {\rm i} \, 1.326~799~292$ & 
$-1.113~340~813~708 \, + \, {\rm i} \, 1.326~827~863~160$ \\      
$10$ & $-11.9~767 \, + \, {\rm i} \, 0.75~063$ & 
$-1.113~334~783 \, + \, {\rm i} \, 1.326~827~026$ & 
$-1.113~340~793~311 \, + \, {\rm i} \, 1.326~827~897~184$ \\      
$11$ & $-6.38~818 \, - \, {\rm i} \, 8.92~900$ & 
$-1.113~340~291 \, + \, {\rm i} \, 1.326~829~768$ & 
$-1.113~340~798~965 \, + \, {\rm i} \, 1.326~827~896~881$ \\      
$12$ & $ 5.43~724 \, - \, {\rm i} \, 8.92~900$ & 
$-1.113~341~213 \, + \, {\rm i} \, 1.326~827~964$ 
& $-1.113~340~798~476 \, + \, {\rm i} \, 1.326~827~896~233$ \\      
$13$ & $ 11.6~655 \, + \, {\rm i} \, 1.85~859$ & 
$-1.113~340~832 \, + \, {\rm i} \, 1.326~827~770$ & 
$-1.113~340~798~438 \, + \, {\rm i} \, 1.326~827~896~343$ \\      
$14$ & $ 5.12~954 \, + \, {\rm i} \, 13.1~791$ & 
$-1.113~340~767 \, + \, {\rm i} \, 1.326~827~891$ & 
$-1.113~340~798~455 \, + \, {\rm i} \, 1.326~827~896~339$ \\      
$15$ & $-8.54~354 \, + \, {\rm i} \, 13.1~791$ & 
$-1.113~340~796 \, + \, {\rm i} \, 1.326~827~905$ & 
$-1.113~340~798~453 \, + \, {\rm i} \, 1.326~827~896~338$ \\       
\hline
Exact & & $- 1.113~340~798 \, + \, {\rm i} \, 1.326~827~896$
& $- 1.113~340~798~453 \, + \, {\rm i} \, 1.326~827~896~338$%
\rule[-1pt]{0pt}{12pt} \\
\end{tabular}
\end{table}

\begin{table}
\caption{Acceleration of the convergence of the hypergeometric series 
${}_2 F_1 (3/7, 5/2; 7/2; z)$ for $z = 77/100$ with the help of Wynn's
epsilon algorithm, Brezinski's theta algorithm, and the Levin-type
transformations $d_{k}^{(n)} (\zeta, s_n)$ and $\delta_{k}^{(n)} (\zeta,
s_n)$.}
\label{Tab_4_1}
\begin{tabular}{lrrrrr}%
$n$%
& \multicolumn{1}{c}{$s_{n} (z)$}%
& \multicolumn{1}{c}{$\epsilon_{2 \Ent {n/2}}^{(n - 2 \Ent {n/2})}$}%
& \multicolumn{1}{c}{$\theta_{2 \Ent {n/3}}^{(n - 3 \Ent {n/3})}$}%
& \multicolumn{1}{c}{$d_{n}^{(0)} \bigl(1, s_0 (z) \bigr)$}%
& \multicolumn{1}{c}{${\delta}_{n}^{(0)} \bigl(1, s_0 (z) \bigr)$}%
\rule[-6pt]{0pt}{18pt} \\
& \multicolumn{1}{c}{Eq.\ (\protect\ref{ParSumT4_1})}%
& \multicolumn{1}{c}{Eq.\ (\protect\ref{eps_al})}%
& \multicolumn{1}{c}{Eq.\ (\protect\ref{theta_al})}%
& \multicolumn{1}{c}{Eq.\ (\protect\ref{dLevTr})}%
& \multicolumn{1}{c}{Eq.\ (\protect\ref{dSidTr})}%
\rule[-6pt]{0pt}{12pt} 
\\
\hline%
$0 $ & $1.000~000~000~000$ & $1.000~000~000~000$ & $1.000~000~000~000$ 
& $1.000~000~000~000$ & $1.000~000~000~000$%
\rule[-1pt]{0pt}{12pt} \\
$1 $ & $1.235~714~285~714$ & $1.235~714~285~714$ & $1.235~714~285~714$ 
& $1.411~927~877~947$ & $1.411~927~877~947$ \\
$2 $ & $1.336~547~619~048$ & $1.411~927~877~947$ & $1.336~547~619~048$ 
& $1.448~820~807~146$ & $1.448~820~807~146$ \\
$3 $ & $1.387~972~619~048$ & $1.441~496~598~639$ & $1.485~041~451~400$ 
& $1.459~648~832~810$ & $1.460~767~058~064$ \\
$4 $ & $1.416~691~503~663$ & $1.457~490~562~327$ & $1.471~789~761~447$ 
& $1.462~687~260~952$ & $1.463~271~135~352$ \\
$5 $ & $1.433~666~279~063$ & $1.461~253~523~425$ & $1.467~349~021~834$ 
& $1.463~515~006~731$ & $1.463~725~027~496$ \\
$6 $ & $1.444~100~773~353$ & $1.463~029~920~830$ & $1.463~655~490~309$ 
& $1.463~733~041~236$ & $1.463~795~467~606$ \\
$7 $ & $1.450~702~748~500$ & $1.463~503~143~435$ & $1.463~758~479~471$ 
& $1.463~788~732~403$ & $1.463~805~515~965$ \\
$8 $ & $1.454~973~597~660$ & $1.463~711~233~069$ & $1.463~790~047~799$ 
& $1.463~802~617~705$ & $1.463~806~889~463$ \\
$9 $ & $1.457~785~549~506$ & $1.463~770~364~632$ & $1.463~806~019~078$ 
& $1.463~806~018~807$ & $1.463~807~072~295$ \\
$10$ & $1.459~663~708~382$ & $1.463~795~262~347$ & $1.463~806~921~084$ 
& $1.463~806~841~297$ & $1.463~807~096~230$ \\
$11$ & $1.460~933~204~660$ & $1.463~802~626~156$ & $1.463~807~092~160$ 
& $1.463~807~038~364$ & $1.463~807~099~326$ \\
$12$ & $1.461~799~964~188$ & $1.463~805~637~323$ & $1.463~807~125~890$ 
& $1.463~807~085~259$ & $1.463~807~099~724$ \\
$13$ & $1.462~396~866~847$ & $1.463~806~552~580$ & $1.463~807~118~461$ 
& $1.463~807~096~360$ & $1.463~807~099~774$ \\
$14$ & $1.462~811~003~597$ & $1.463~806~919~041$ & $1.463~807~119~752$ 
& $1.463~807~098~978$ & $1.463~807~099~781$ \\
$15$ & $1.463~100~213~038$ & $1.463~807~032~663$ & $1.463~807~119~214$ 
& $1.463~807~099~593$ & $1.463~807~099~781$ \\
$16$ & $1.463~303~343~588$ & $1.463~807~077~439$ & $1.463~807~118~928$ 
& $1.463~807~099~737$ & $1.463~807~099~782$ \\
\hline
Exact & $1.463~807~099~782$ & $1.463~807~099~782$ & $1.463~807~099~782$ 
& $1.463~807~099~782$ & $1.463~807~099~782$%
\rule[-1pt]{0pt}{12pt} \\
\end{tabular}
\end{table}

\begin{table}
\caption{Acceleration of the convergence of the hypergeometric series 
${}_2 F_1 (3/7, 5/2; - 7/2; z)$ for $z = 77/100$ with the help of 
Aitken's iterated $\Delta^2$ process, Wynn's epsilon algorithm, and 
Brezinski's theta algorithm.}
\label{Tab_4_2}
\begin{tabular}{lrrrr}%
$n$%
& \multicolumn{1}{c}{$s_{n} (z)$}%
& \multicolumn{1}{c}{${\cal A}_{\Ent {n/2}}^{(n - 2 \Ent {n/2})}$}%
& \multicolumn{1}{c}{$\epsilon_{2 \Ent {n/2}}^{(n - 2 \Ent {n/2})}$}%
& \multicolumn{1}{c}{$\theta_{2 \Ent {n/3}}^{(n - 3 \Ent {n/3})}$}%
\rule[-6pt]{0pt}{18pt} \\
& \multicolumn{1}{c}{Eq.\ (\protect\ref{ParSumT4_2})}%
& \multicolumn{1}{c}{Eq.\ (\protect\ref{It_Aitken})}%
& \multicolumn{1}{c}{Eq.\ (\protect\ref{eps_al})}%
& \multicolumn{1}{c}{Eq.\ (\protect\ref{theta_al})}%
\rule[-6pt]{0pt}{12pt} 
\\
\hline%
$0 $ & $ 1.000~000 \cdot 10^{+0}$ & $ 1.000~000~000~000 \cdot 10^{+0}$ 
& $ 1.000~000~000~000 \cdot 10^{+0}$ & $ 1.000~000~000~000 \cdot 10^{+0}$%
\rule[-1pt]{0pt}{12pt} \\
$1 $ & $ 7.642~857 \cdot 10^{-1}$ & $ 7.642~857~142~857 \cdot 10^{-1}$ 
& $ 7.642~857~142~857 \cdot 10^{-1}$ & $ 7.642~857~142~857 \cdot 10^{-1}$ \\
$2 $ & $ 9.457~857 \cdot 10^{-1}$ & $ 8.668~280~871~671 \cdot 10^{-1}$ 
& $ 8.668~280~871~671 \cdot 10^{-1}$ & $ 9.457~857~142~857 \cdot 10^{-1}$ \\
$3 $ & $ 6.063~807 \cdot 10^{-1}$ & $ 8.275~261~324~042 \cdot 10^{-1}$ 
& $ 8.275~261~324~042 \cdot 10^{-1}$ & $ 8.485~759~131~562 \cdot 10^{-1}$ \\
$4 $ & $ 3.070~461 \cdot 10^{+0}$ & $ 8.535~660~811~732 \cdot 10^{-1}$ 
& $ 8.609~189~507~196 \cdot 10^{-1}$ & $ 8.494~334~311~697 \cdot 10^{-1}$ \\
$5 $ & $ 2.491~700 \cdot 10^{+1}$ & $ 8.361~725~059~741 \cdot 10^{-1}$ 
& $ 8.157~392~107~013 \cdot 10^{-1}$ & $ 8.453~777~105~530 \cdot 10^{-1}$ \\
$6 $ & $ 1.010~158 \cdot 10^{+2}$ & $ 8.516~151~469~394 \cdot 10^{-1}$ 
& $ 8.896~493~143~013 \cdot 10^{-1}$ & $ 8.487~841~750~682 \cdot 10^{-1}$ \\
$7 $ & $ 2.839~789 \cdot 10^{+2}$ & $ 8.209~992~692~837 \cdot 10^{-1}$ 
& $ 4.152~615~385~984 \cdot 10^{-1}$ & $ 8.486~139~776~736 \cdot 10^{-1}$ \\
$8 $ & $ 6.390~582 \cdot 10^{+2}$ & $ 8.549~782~290~518 \cdot 10^{-1}$ 
& $ 1.979~200~475~641 \cdot 10^{+0}$ & $ 8.486~753~544~589 \cdot 10^{-1}$ \\
$9 $ & $ 1.236~512 \cdot 10^{+3}$ & $ 9.170~319~462~988 \cdot 10^{-1}$ 
& $ 1.579~725~923~915 \cdot 10^{+1}$ & $ 8.486~400~009~032 \cdot 10^{-1}$ \\
$10$ & $ 2.143~447 \cdot 10^{+3}$ & $ 8.324~631~947~362 \cdot 10^{-1}$ 
& $-9.032~991~315~812 \cdot 10^{+1}$ & $ 8.486~226~769~041 \cdot 10^{-1}$ \\
$11$ & $ 3.416~644 \cdot 10^{+3}$ & $ 7.319~750~281~906 \cdot 10^{-1}$ 
& $-4.871~173~637~209 \cdot 10^{+2}$ & $ 8.486~314~552~725 \cdot 10^{-1}$ \\
$12$ & $ 5.097~263 \cdot 10^{+3}$ & $ 8.296~736~361~657 \cdot 10^{-1}$ 
& $ 9.550~458~201~456 \cdot 10^{+3}$ & $ 8.486~265~235~534 \cdot 10^{-1}$ \\
$13$ & $ 7.207~772 \cdot 10^{+3}$ & $ 2.786~711~790~541 \cdot 10^{+0}$ 
& $ 3.152~031~922~861 \cdot 10^{+4}$ & $ 8.486~238~887~649 \cdot 10^{-1}$ \\
$14$ & $ 9.751~014 \cdot 10^{+3}$ & $ 2.417~040~069~989 \cdot 10^{+0}$ 
& $ 9.017~022~343~102 \cdot 10^{+4}$ & $ 8.486~245~607~105 \cdot 10^{-1}$ \\
$15$ & $ 1.271~111 \cdot 10^{+4}$ & $ 2.660~054~094~715 \cdot 10^{+0}$ 
& $ 9.807~349~664~879 \cdot 10^{+4}$ & $ 8.486~241~968~013 \cdot 10^{-1}$ \\
$16$ & $ 1.605~569 \cdot 10^{+4}$ & $ 1.583~428~788~102 \cdot 10^{+0}$ 
& $ 1.006~124~741~024 \cdot 10^{+5}$ & $ 8.486~240~157~205 \cdot 10^{-1}$ \\
$17$ & $ 1.973~906 \cdot 10^{+4}$ & $ 2.773~669~575~356 \cdot 10^{+0}$ 
& $ 1.009~114~496~621 \cdot 10^{+5}$ & $ 8.485~741~255~417 \cdot 10^{-1}$ \\
$18$ & $ 2.370~573 \cdot 10^{+4}$ & $ 2.553~560~837~044 \cdot 10^{+0}$ 
& $ 1.009~964~826~770 \cdot 10^{+5}$ & $ 8.486~242~524~145 \cdot 10^{-1}$ \\
$19$ & $ 2.789~405 \cdot 10^{+4}$ & $ 2.384~039~609~516 \cdot 10^{+0}$ 
& $ 1.010~111~063~638 \cdot 10^{+5}$ & $ 8.486~355~858~960 \cdot 10^{-1}$ \\
$20$ & $ 3.223~965 \cdot 10^{+4}$ & $ 2.454~155~706~579 \cdot 10^{+0}$ 
& $ 1.010~153~078~626 \cdot 10^{+5}$ & $ 8.492~820~131~138 \cdot 10^{-1}$ \\
$21$ & $ 3.667~836 \cdot 10^{+4}$ & $ 2.384~021~379~340 \cdot 10^{+0}$ 
& $ 1.010~161~882~819 \cdot 10^{+5}$ & $ 8.486~069~295~394 \cdot 10^{-1}$ \\
$22$ & $ 4.114~878 \cdot 10^{+4}$ & $ 2.454~149~413~545 \cdot 10^{+0}$ 
& $ 1.010~164~516~908 \cdot 10^{+5}$ & $ 8.470~533~544~441 \cdot 10^{-1}$ \\
$23$ & $ 4.559~414 \cdot 10^{+4}$ & $ 7.889~100~799~808 \cdot 10^{+2}$ 
& $ 1.010~165~130~151 \cdot 10^{+5}$ & $ 2.015~803~318~335 \cdot 10^{+0}$ \\
$24$ & $ 4.996~370 \cdot 10^{+4}$ & $ 7.855~297~113~412 \cdot 10^{+2}$ 
& $ 1.010~165~321~636 \cdot 10^{+5}$ & $ 8.519~514~381~837 \cdot 10^{-1}$ \\
$25$ & $ 5.421~360 \cdot 10^{+4}$ & $ 7.916~396~789~254 \cdot 10^{+2}$ 
& $ 1.010~165~369~332 \cdot 10^{+5}$ & $ 1.185~791~071~696 \cdot 10^{+0}$ \\
$26$ & $ 5.830~726 \cdot 10^{+4}$ & $ 7.874~382~737~052 \cdot 10^{+2}$ 
& $ 1.010~165~384~712 \cdot 10^{+5}$ & $-6.192~361~264~079 \cdot 10^{+0}$ \\
$27$ & $ 6.221~545 \cdot 10^{+4}$ & $ 7.780~517~151~440 \cdot 10^{+2}$ 
& $ 1.010~165~388~746 \cdot 10^{+5}$ & $-2.153~607~969~121 \cdot 10^{+0}$ \\
$28$ & $ 6.591~599 \cdot 10^{+4}$ & $ 7.781~235~296~899 \cdot 10^{+2}$ 
& $ 1.010~165~390~075 \cdot 10^{+5}$ & $-7.845~642~420~775 \cdot 10^{+0}$ \\
$29$ & $ 6.939~333 \cdot 10^{+4}$ & $ 7.780~882~446~428 \cdot 10^{+2}$ 
& $ 1.010~165~390~438 \cdot 10^{+5}$ & $ 2.998~151~673~644 \cdot 10^{+2}$ \\
$30$ & $ 7.263~789 \cdot 10^{+4}$ & $ 7.781~060~228~285 \cdot 10^{+2}$ 
& $ 1.010~165~390~560 \cdot 10^{+5}$ & $ 3.199~681~242~234 \cdot 10^{+1}$ \\
\hline
Exact & $1.010~165 \cdot 10^{+5}$ & $1.010~165~390~611 \cdot 10^{+5}$ 
& $1.010~165~390~611 \cdot 10^{+5}$ 
& $1.010~165~390~611 \cdot 10^{+5}$\rule[-1pt]{0pt}{12pt} \\
\end{tabular}
\end{table}

\begin{table}
\caption{Acceleration of the convergence of the hypergeometric series 
${}_2 F_1 (3/7, 5/2; - 7/2; z)$ for $z = 77/100$ with the help of the
iteration of Brezinski's theta algorithm and the Levin-type
transformations $d_{k}^{(n)} (\zeta, s_n)$ and $\delta_{k}^{(n)} (\zeta,
s_n)$.}
\label{Tab_4_2_a}
\begin{tabular}{lrrrr}%
$n$%
& \multicolumn{1}{c}{$s_{n} (z)$}%
& \multicolumn{1}{c}{${\cal J}_{\Ent {n/3}}^{(n - 3 \Ent {n/3})}$}%
& \multicolumn{1}{c}{$d_{n}^{(0)} \bigl(1, s_0 (z) \bigr)$}%
& \multicolumn{1}{c}{${\delta}_{n}^{(0)} \bigl(1, s_0 (z) \bigr)$}%
\rule[-6pt]{0pt}{18pt} \\
& \multicolumn{1}{c}{Eq.\ (\protect\ref{ParSumT4_2})}%
& \multicolumn{1}{c}{Eq.\ (\protect\ref{thetit_1})}%
& \multicolumn{1}{c}{Eq.\ (\protect\ref{dLevTr})}%
& \multicolumn{1}{c}{Eq.\ (\protect\ref{dSidTr})}%
\rule[-6pt]{0pt}{12pt} 
\\
\hline%
$0 $ & $ 1.000~000 \cdot 10^{+0}$ & $ 1.000~000~000~000 \cdot 10^{+0}$ &
$ 1.000~000~000~000 \cdot 10^{+0}$ & $ 1.000~000~000~000 \cdot 10^{+0}$%
\rule[-1pt]{0pt}{12pt} \\
$1 $ & $ 7.642~857 \cdot 10^{-1}$ & $ 7.642~857~142~857 \cdot 10^{-1}$ 
& $ 8.668~280~871~671 \cdot 10^{-1}$ & $ 8.668~280~871~671 \cdot 10^{-1}$ \\
$2 $ & $ 9.457~857 \cdot 10^{-1}$ & $ 9.457~857~142~857 \cdot 10^{-1}$ 
& $ 8.384~394~404~076 \cdot 10^{-1}$ & $ 8.384~394~404~076 \cdot 10^{-1}$ \\
$3 $ & $ 6.063~807 \cdot 10^{-1}$ & $ 8.485~759~131~562 \cdot 10^{-1}$ 
& $ 8.564~401~754~536 \cdot 10^{-1}$ & $ 8.543~029~016~450 \cdot 10^{-1}$ \\
$4 $ & $ 3.070~461 \cdot 10^{+0}$ & $ 8.494~334~311~697 \cdot 10^{-1}$ 
& $ 8.410~428~387~954 \cdot 10^{-1}$ & $ 8.451~555~013~928 \cdot 10^{-1}$ \\
$5 $ & $ 2.491~700 \cdot 10^{+1}$ & $ 8.453~777~105~530 \cdot 10^{-1}$ 
& $ 8.581~719~847~922 \cdot 10^{-1}$ & $ 8.511~715~495~996 \cdot 10^{-1}$ \\
$6 $ & $ 1.010~158 \cdot 10^{+2}$ & $ 8.487~877~036~926 \cdot 10^{-1}$ 
& $ 8.339~732~531~530 \cdot 10^{-1}$ & $ 8.466~755~128~315 \cdot 10^{-1}$ \\
$7 $ & $ 2.839~789 \cdot 10^{+2}$ & $ 8.485~959~518~738 \cdot 10^{-1}$ 
& $ 8.765~035~290~451 \cdot 10^{-1}$ & $ 8.504~596~561~867 \cdot 10^{-1}$ \\
$8 $ & $ 6.390~582 \cdot 10^{+2}$ & $ 8.489~925~522~131 \cdot 10^{-1}$ 
& $ 7.853~488~289~879 \cdot 10^{-1}$ & $ 8.469~220~701~568 \cdot 10^{-1}$ \\
$9 $ & $ 1.236~512 \cdot 10^{+3}$ & $ 8.486~664~090~804 \cdot 10^{-1}$ 
& $ 1.019~118~632~891 \cdot 10^{+0}$ & $ 8.505~456~725~376 \cdot 10^{-1}$ \\
$10$ & $ 2.143~447 \cdot 10^{+3}$ & $ 8.335~704~458~916 \cdot 10^{-1}$ 
& $ 3.141~144~923~964 \cdot 10^{-1}$ & $ 8.465~273~890~370 \cdot 10^{-1}$ \\
$11$ & $ 3.416~644 \cdot 10^{+3}$ & $ 8.497~498~135~494 \cdot 10^{-1}$ 
& $ 2.776~012~986~506 \cdot 10^{+0}$ & $ 8.513~020~606~999 \cdot 10^{-1}$ \\
$12$ & $ 5.097~263 \cdot 10^{+3}$ & $ 8.501~342~364~535 \cdot 10^{-1}$ 
& $-7.045~936~018~531 \cdot 10^{+0}$ & $ 8.452~752~921~910 \cdot 10^{-1}$ \\
$13$ & $ 7.207~772 \cdot 10^{+3}$ & $ 8.498~411~613~838 \cdot 10^{-1}$ 
& $ 3.717~678~472~983 \cdot 10^{+1}$ & $ 8.532~975~348~760 \cdot 10^{-1}$ \\
$14$ & $ 9.751~014 \cdot 10^{+3}$ & $ 8.497~632~848~220 \cdot 10^{-1}$ 
& $-1.858~061~858~604 \cdot 10^{+2}$ & $ 8.421~055~650~264 \cdot 10^{-1}$ \\
$15$ & $ 1.271~111 \cdot 10^{+4}$ & $ 8.497~961~772~597 \cdot 10^{-1}$ 
& $ 1.044~967~163~943 \cdot 10^{+3}$ & $ 8.583~855~953~130 \cdot 10^{-1}$ \\
$16$ & $ 1.605~569 \cdot 10^{+4}$ & $ 8.498~474~972~833 \cdot 10^{-1}$ 
& $-7.005~065~553~184 \cdot 10^{+3}$ & $ 8.338~036~004~332 \cdot 10^{-1}$ \\
$17$ & $ 1.973~906 \cdot 10^{+4}$ & $ 7.503~910~591~903 \cdot 10^{-1}$ 
& $ 3.081~110~326~042 \cdot 10^{+4}$ & $ 8.721~871~069~015 \cdot 10^{-1}$ \\
$18$ & $ 2.370~573 \cdot 10^{+4}$ & $ 8.499~855~007~426 \cdot 10^{-1}$ 
& $ 1.465~687~756~463 \cdot 10^{+5}$ & $ 8.104~110~955~846 \cdot 10^{-1}$ \\
$19$ & $ 2.789~405 \cdot 10^{+4}$ & $ 7.414~199~029~835 \cdot 10^{-1}$ 
& $ 9.719~482~708~581 \cdot 10^{+4}$ & $ 9.126~025~807~834 \cdot 10^{-1}$ \\
$20$ & $ 3.223~965 \cdot 10^{+4}$ & $ 6.852~823~777~347 \cdot 10^{-1}$ 
& $ 1.014~869~541~810 \cdot 10^{+5}$ & $ 7.392~774~724~769 \cdot 10^{-1}$ \\
$21$ & $ 3.667~836 \cdot 10^{+4}$ & $ 6.800~566~895~738 \cdot 10^{-1}$ 
& $ 1.009~649~257~196 \cdot 10^{+5}$ & $ 1.040~040~245~220 \cdot 10^{+0}$ \\
$22$ & $ 4.114~878 \cdot 10^{+4}$ & $ 6.778~745~199~478 \cdot 10^{-1}$ 
& $ 1.010~218~868~540 \cdot 10^{+5}$ & $ 5.070~943~906~098 \cdot 10^{-1}$ \\
$23$ & $ 4.559~414 \cdot 10^{+4}$ & $ 6.776~231~440~878 \cdot 10^{-1}$ 
& $ 1.010~160~173~901 \cdot 10^{+5}$ & $ 1.469~832~601~187 \cdot 10^{+0}$ \\
$24$ & $ 4.996~370 \cdot 10^{+4}$ & $ 6.777~792~652~237 \cdot 10^{-1}$ 
& $ 1.010~165~871~479 \cdot 10^{+5}$ & $-3.005~355~825~647 \cdot 10^{-1}$ \\
$25$ & $ 5.421~360 \cdot 10^{+4}$ & $ 6.807~217~929~960 \cdot 10^{-1}$ 
& $ 1.010~165~348~613 \cdot 10^{+5}$ & $ 3.008~947~036~607 \cdot 10^{+0}$ \\
$26$ & $ 5.830~726 \cdot 10^{+4}$ & $ 6.775~201~071~781 \cdot 10^{-1}$ 
& $ 1.010~165~394~094 \cdot 10^{+5}$ & $-3.273~286~658~302 \cdot 10^{+0}$ \\
$27$ & $ 6.221~545 \cdot 10^{+4}$ & $ 6.796~867~331~534 \cdot 10^{-1}$ 
& $ 1.010~165~390~335 \cdot 10^{+5}$ & $ 8.821~765~469~703 \cdot 10^{+0}$ \\
$28$ & $ 6.591~599 \cdot 10^{+4}$ & $ 1.050~749~944~450 \cdot 10^{+0}$ 
& $ 1.010~165~390~631 \cdot 10^{+5}$ & $-1.477~966~683~169 \cdot 10^{+1}$ \\
$29$ & $ 6.939~333 \cdot 10^{+4}$ & $-6.085~082~692~546 \cdot 10^{-1}$ 
& $ 1.010~165~390~609 \cdot 10^{+5}$ & $ 3.183~343~913~267 \cdot 10^{+1}$ \\
$30$ & $ 7.263~789 \cdot 10^{+4}$ & $-6.037~663~167~562 \cdot 10^{-1}$ 
& $ 1.010~165~390~611 \cdot 10^{+5}$ & $-6.136~838~085~955 \cdot 10^{+1}$ \\
\hline
Exact & $1.010~165 \cdot 10^{+5}$ & $1.010~165~390~611 \cdot 10^{+5}$ 
& $1.010~165~390~611 \cdot 10^{+5}$ 
& $1.010~165~390~611 \cdot 10^{+5}$\rule[-1pt]{0pt}{12pt} \\
\end{tabular}
\end{table}

\begin{table}
\caption{Acceleration of the convergence of the hypergeometric series 
${}_2 F_1 (3/7, 5/2; - 7/2; z)$ for $z = 77/100$ with the help of Wynn's
epsilon algorithm and the Levin-type transformations $d_{k}^{(n)}
(\zeta, s_n)$ and $\delta_{k}^{(n)} (\zeta, s_n)$.}
\label{Tab_4_3}
\begin{tabular}{lrrrr}%
$n$%
& \multicolumn{1}{c}{$s_{n}^{(22)} (z)$}%
& \multicolumn{1}{c}{$\epsilon_{2 \Ent {n/2}}^{(n - 2 \Ent {n/2})}$}%
& \multicolumn{1}{c}{$d_{n}^{(0)} \bigl(1, s_{0}^{(22)} (z) \bigr)$}%
& \multicolumn{1}{c}
{${\delta}_{n}^{(0)} \bigl(1, s_{0}^{(22)} (z) \bigr)$}%
\rule[-6pt]{0pt}{18pt} \\
& \multicolumn{1}{c}{Eq.\ (\protect\ref{ParSumT4_3})}%
& \multicolumn{1}{c}{Eq.\ (\protect\ref{eps_al})}%
& \multicolumn{1}{c}{Eq.\ (\protect\ref{dLevTr})}%
& \multicolumn{1}{c}{Eq.\ (\protect\ref{dSidTr})}%
\rule[-6pt]{0pt}{12pt} 
\\
\hline%
$0 $ & $4.114~878 \cdot 10^{4}$ & $4.114~877~834~479 \cdot 10^{4}$ & 
$4.114~877~834~479 \cdot 10^{4}$ & $4.114~877~834~479 \cdot 10^{4}$%
\rule[-1pt]{0pt}{12pt} \\
$1 $ & $4.559~414 \cdot 10^{4}$	& $4.559~413~879~729 \cdot 10^{4}$ & 
$3.018~541~281~910 \cdot 10^{5}$ & $3.018~541~281~910 \cdot 10^{5}$ \\
$2 $ & $4.996~370 \cdot 10^{4}$	& $3.018~541~281~910 \cdot 10^{5}$ & 
$1.806~947~255~600 \cdot 10^{5}$ & $1.806~947~255~600 \cdot 10^{5}$ \\
$3 $ & $5.421~360 \cdot 10^{4}$	& $2.051~523~999~854 \cdot 10^{5}$ & 
$1.359~244~863~515 \cdot 10^{5}$ & $1.297~795~854~790 \cdot 10^{5}$ \\
$4 $ & $5.830~726 \cdot 10^{4}$	& $8.720~589~218~810 \cdot 10^{4}$ & 
$1.163~496~269~997 \cdot 10^{5}$ & $1.101~826~174~021 \cdot 10^{5}$ \\
$5 $ & $6.221~545 \cdot 10^{4}$	& $9.074~241~989~747 \cdot 10^{4}$ & 
$1.074~151~889~285 \cdot 10^{5}$ & $1.032~936~136~116 \cdot 10^{5}$ \\
$6 $ & $6.591~599 \cdot 10^{4}$	& $1.032~456~488~505 \cdot 10^{5}$ & 
$1.034~564~538~627 \cdot 10^{5}$ & $1.013~474~447~563 \cdot 10^{5}$ \\
$7 $ & $6.939~333 \cdot 10^{4}$	& $1.025~201~938~923 \cdot 10^{5}$ & 
$1.018~351~267~311 \cdot 10^{5}$ & $1.009~927~590~147 \cdot 10^{5}$ \\
$8 $ & $7.263~789 \cdot 10^{4}$	& $1.008~495~311~117 \cdot 10^{5}$ & 
$1.012~445~529~468 \cdot 10^{5}$ & $1.009~872~220~781 \cdot 10^{5}$ \\
$9 $ & $7.564~538 \cdot 10^{4}$	& $1.009~077~057~226 \cdot 10^{5}$ & 
$1.010~614~403~464 \cdot 10^{5}$ & $1.010~079~938~133 \cdot 10^{5}$ \\
$10$ & $7.841~602 \cdot 10^{4}$	& $1.010~227~941~973 \cdot 10^{5}$ & 
$1.010~172~932~578 \cdot 10^{5}$ & $1.010~153~214~789 \cdot 10^{5}$ \\
$11$ & $8.095~384 \cdot 10^{4}$	& $1.010~204~695~588 \cdot 10^{5}$ & 
$1.010~117~444~530 \cdot 10^{5}$ & $1.010~165~544~386 \cdot 10^{5}$ \\
$12$ & $8.326~588 \cdot 10^{4}$	& $1.010~164~792~364 \cdot 10^{5}$ & 
$1.010~135~304~199 \cdot 10^{5}$ & $1.010~165~871~909 \cdot 10^{5}$ \\
$13$ & $8.536~159 \cdot 10^{4}$	& $1.010~165~026~181 \cdot 10^{5}$ & 
$1.010~152~575~006 \cdot 10^{5}$ & $1.010~165~498~236 \cdot 10^{5}$ \\
$14$ & $8.725~217 \cdot 10^{4}$	& $1.010~165~383~969 \cdot 10^{5}$ & 
$1.010~160~984~017 \cdot 10^{5}$ & $1.010~165~400~397 \cdot 10^{5}$ \\
$15$ & $8.895~003 \cdot 10^{4}$	& $1.010~165~386~674 \cdot 10^{5}$ & 
$1.010~164~122~704 \cdot 10^{5}$ & $1.010~165~389~960 \cdot 10^{5}$ \\
$16$ & $9.046~835 \cdot 10^{4}$	& $1.010~165~390~431 \cdot 10^{5}$ & 
$1.010~165~091~529 \cdot 10^{5}$ & $1.010~165~390~292 \cdot 10^{5}$ \\
$17$ & $9.182~068 \cdot 10^{4}$	& $1.010~165~390~507 \cdot 10^{5}$ & 
$1.010~165~339~419 \cdot 10^{5}$ & $1.010~165~390~569 \cdot 10^{5}$ \\
$18$ & $9.302~057 \cdot 10^{4}$	& $1.010~165~390~603 \cdot 10^{5}$ & 
$1.010~165~388~497 \cdot 10^{5}$ & $1.010~165~390~609 \cdot 10^{5}$ \\
$19$ & $9.408~137 \cdot 10^{4}$	& $1.010~165~390~607 \cdot 10^{5}$ & 
$1.010~165~393~632 \cdot 10^{5}$ & $1.010~165~390~611 \cdot 10^{5}$ \\
$20$ & $9.501~598 \cdot 10^{4}$	& $1.010~165~390~610 \cdot 10^{5}$ & 
$1.010~165~392~321 \cdot 10^{5}$ & $1.010~165~390~611 \cdot 10^{5}$ \\
\hline
Exact & $1.010~165 \cdot 10^{5}$ & $1.010~165~390~611 \cdot 10^{5}$ 
& $1.010~165~390~611 \cdot 10^{5}$ 
& $1.010~165~390~611 \cdot 10^{5}$\rule[-1pt]{0pt}{12pt} \\
\end{tabular}
\end{table}

\narrowtext

\begin{table}
\caption{Selected coefficients $c_{n}^{(3)}$ of the renormalized weak
coupling expansion (\protect\ref{rwcEm}) as well as the corresponding
ratios ${\cal C}_{n}^{(3)}$ defined by Eq.\ (\protect\ref{Asy_n_C3}).}
\label{Tab_b_1}
\begin{tabular}{lr@{.}lr}%
$n$%
& \multicolumn{2}{c}{$c_{n}^{(3)}$}%
& \multicolumn{1}{c}{$\quad {\cal C}_{n}^{(3)}$}%
\rule[-6pt]{0pt}{18pt} 
\\
\hline%
$0  $ & $ 1$ & $000~000~000                 $ & $-0.49218$ \\
$1  $ & $-3$ & $333~333~333 \cdot 10^{-1}   $ & $-1.51799$ \\
$2  $ & $-9$ & $074~074~074 \cdot 10^{-2}   $ & $ 0.32773$ \\
$3  $ & $ 3$ & $451~646~091 \cdot 10^{-1}   $ & $ 0.34954$ \\
$4  $ & $-3$ & $064~808~585                 $ & $ 0.44181$ \\
$5  $ & $ 4$ & $145~321~167 \cdot 10^{+1}   $ & $ 0.51354$ \\
$6  $ & $-8$ & $011~680~849 \cdot 10^{+2}   $ & $ 0.57041$ \\
$7  $ & $ 2$ & $103~995~759 \cdot 10^{+04}  $ & $ 0.61602$ \\
$8  $ & $-7$ & $225~346~394 \cdot 10^{+05}  $ & $ 0.65320$ \\
$9  $ & $ 3$ & $148~105~306 \cdot 10^{+07}  $ & $ 0.68398$ \\
$10 $ & $-1$ & $698~432~299 \cdot 10^{+09}  $ & $ 0.70986$ \\
$11 $ & $ 1$ & $112~192~278 \cdot 10^{+11}  $ & $ 0.73188$ \\
$12 $ & $-8$ & $693~791~326 \cdot 10^{+12}  $ & $ 0.75085$ \\
$13 $ & $ 7$ & $998~709~458 \cdot 10^{+14}  $ & $ 0.76734$ \\
$14 $ & $-8$ & $558~133~512 \cdot 10^{+16}  $ & $ 0.78181$ \\
$15 $ & $ 1$ & $053~809~185 \cdot 10^{+19}  $ & $ 0.79460$ \\
$25 $ & $ 3$ & $966~243~637 \cdot 10^{+42}  $ & $ 0.87072$ \\
$50 $ & $-8$ & $551~437~639 \cdot 10^{+114} $ & $ 0.93301$ \\
$75 $ & $ 4$ & $053~097~428 \cdot 10^{+198} $ & $ 0.95481$ \\
$100$ & $-4$ & $540~614~326 \cdot 10^{+289} $ & $ 0.96591$ \\
$125$ & $ 1$ & $552~803~192 \cdot 10^{+386} $ & $ 0.97263$ \\
$150$ & $-1$ & $248~762~212 \cdot 10^{+487} $ & $ 0.97714$ \\
$175$ & $ 4$ & $338~412~567 \cdot 10^{+591} $ & $ 0.98037$ \\
$200$ & $-1$ & $952~374~463 \cdot 10^{+699} $ & $ 0.98280$ \\
$225$ & $ 4$ & $623~901~073 \cdot 10^{+809} $ & $ 0.98470$ \\
$250$ & $-2$ & $864~322~271 \cdot 10^{+922} $ & $ 0.98622$ \\
$275$ & $ 2$ & $654~720~800 \cdot 10^{+1037}$ & $ 0.98746$ \\
$300$ & $-2$ & $331~943~009 \cdot 10^{+1154}$ & $ 0.98850$ \\
\end{tabular}
\end{table}

\begin{table}
\caption{Infinite coupling limit $k_3$ of the sextic anharmonic
oscillator. Convergence of the approximants $k_{3}^{(n, l)}$ defined by
Eq.\ (\protect\ref{k3_ln}) with the three highest possible values of $n
\le 299 - l$ for $0 \le l \le 12$.}
\label{Tab_b_2}
\begin{tabular}{l|rr}%
$l$%
& \multicolumn{1}{c}{$n$}%
& \multicolumn{1}{c}{$k_{3}^{(n, l)}$}%
\\
\hline%
$0$  & $297$ & $1.144~802~453~797~051~992~831$ \\	 
     & $298$ & $1.144~802~453~797~052~042~002$ \\	 
     & $299$ & $1.144~802~453~797~052~088~337$ \\
\hline              
$1$  & $296$ & $1.144~802~453~797~052~586~551$ \\   
     & $297$ & $1.144~802~453~797~052~601~411$ \\   
     & $298$ & $1.144~802~453~797~052~615~225$ \\
\hline              
$2$  & $295$ & $1.144~802~453~797~052~759~384$ \\   
     & $296$ & $1.144~802~453~797~052~761~860$ \\   
     & $297$ & $1.144~802~453~797~052~764~077$ \\
\hline              
$3$  & $294$ & $1.144~802~453~797~052~783~941$ \\   
     & $295$ & $1.144~802~453~797~052~783~501$ \\   
     & $296$ & $1.144~802~453~797~052~783~043$ \\
\hline              
$4$  & $293$ & $1.144~802~453~797~052~775~728$ \\   
     & $294$ & $1.144~802~453~797~052~775~182$ \\   
     & $295$ & $1.144~802~453~797~052~774~657$ \\
\hline              
$5$  & $292$ & $1.144~802~453~797~052~767~958$ \\   
     & $293$ & $1.144~802~453~797~052~767~714$ \\   
     & $294$ & $1.144~802~453~797~052~767~482$ \\
\hline              
$6$  & $291$ & $1.144~802~453~797~052~764~644$ \\   
     & $292$ & $1.144~802~453~797~052~764~570$ \\   
     & $293$ & $1.144~802~453~797~052~764~501$ \\
\hline              
$7$  & $290$ & $1.144~802~453~797~052~763~692$ \\   
     & $291$ & $1.144~802~453~797~052~763~725$ \\   
     & $292$ & $1.144~802~453~797~052~763~758$ \\
\hline              
$8$  & $289$ & $1.144~802~453~797~052~765~158$ \\   
     & $290$ & $1.144~802~453~797~052~764~915$ \\   
     & $291$ & $1.144~802~453~797~052~764~334$ \\
\hline              
$9$  & $288$ & $1.144~802~453~797~052~746~941$ \\   
     & $289$ & $1.144~802~453~797~052~719~382$ \\   
     & $290$ & $1.144~802~453~797~052~727~384$ \\
\hline	
$10$ & $287$ & $1.144~802~453~797~049~918~777$ \\   
     & $288$ & $1.144~802~453~797~053~441~988$ \\   
     & $289$ & $1.144~802~453~797~055~805~549$ \\
\hline	
$11$ & $286$ & $1.144~802~453~797~422~201~710$ \\   
     & $287$ & $1.144~802~453~797~325~472~977$ \\   
     & $288$ & $1.144~802~453~796~999~625~608$ \\
\hline	
$12$ & $285$ & $1.144~802~453~789~565~166~506$ \\   
     & $286$ & $1.144~802~453~760~484~424~039$ \\   
     & $287$ & $1.144~802~453~769~883~065~332$ \\
\end{tabular}
\end{table}

\end{document}